\documentclass[12pt]{extarticle}
\usepackage{extsizes}
\usepackage{amssymb,amsfonts,amsthm}
\usepackage{amsmath}
\usepackage[makeroom]{cancel}
\usepackage{mathtools}
\usepackage{mathrsfs,pifont}
\usepackage{slashed,mathabx} 
\usepackage[bbgreekl]{mathbbol}

\usepackage[top=1in, bottom=1.25in, left=0.75in,
right=0.75in]{geometry}

\usepackage[font={small,sl}]{caption}
\usepackage[all]{xy}

\newcommand{\C}{\mathbb{C}}

\newcommand{\Z}{\mathbb{Z}}
\newcommand{\Q}{\mathbb{Q}}
\newcommand{\R}{\mathbb{R}}
\newcommand{\N}{\mathbb{N}}

\newcommand{\bB}{\mathsf{B}}
\newcommand{\bY}{\mathsf{Y}}

\newcommand{\bv}{\mathsf{v}}
\newcommand{\bw}{\mathsf{w}}

\newcommand{\bX}{\mathsf{X}}
\newcommand{\bJ}{\mathsf{J}}
\newcommand{\bT}{\mathsf{T}}

\newcommand{\bA}{\mathsf{A}}

\newcommand{\bP}{\mathbb{P}}

\newcommand{\cA}{\mathscr{A}}
\newcommand{\cX}{\mathscr{X}}
\newcommand{\cP}{\mathscr{P}}
\newcommand{\cK}{\mathscr{K}}
\newcommand{\cL}{\mathscr{L}}

\newcommand{\cO}{\mathscr{O}}

\newcommand{\cU}{\mathscr{U}}

\newcommand{\cM}{\mathscr{M}}
\newcommand{\cV}{\mathscr{V}}

\newcommand{\cC}{\mathscr{C}}

\newcommand{\bfM}{\mathbf{M}}

\newcommand{\bfG}{\mathbf{G}}

\newcommand{\fg}{\mathfrak{g}} 
\newcommand{\fgh}{\widehat{\mathfrak{g}}} 
 
\newcommand{\fh}{\mathfrak{h}} 
\newcommand{\fhh}{\widehat{\mathfrak{h}}}

\newcommand{\rd}{/\!\!/\!\!/\!\!/}
\newcommand{\rdd}{/\!\!/}
\newcommand{\Ct}{\mathbb{C}^\times}

\newcommand{\fC}{\mathfrak{C}}
\newcommand{\fK}{\mathfrak{K}}

\newcommand{\Hd}{{H}^{\raisebox{0.5mm}{$\scriptscriptstyle \bullet$}}}

\newcommand{\glh}{\widehat{\mathfrak{gl}}} 
\newcommand{\gh}{\widehat{\mathfrak{g}}}

\newcommand{\bk}{\boldsymbol{\kappa}}

\DeclareMathOperator{\chern}{ch}
\DeclareMathOperator{\Def}{Def}

\DeclareMathOperator{\Obs}{Obs}
\DeclareMathOperator{\Lie}{Lie}

\DeclareMathOperator{\Mat}{Mat}

\DeclareMathOperator{\Hilb}{Hilb}
\DeclareMathOperator{\End}{End}

\DeclareMathOperator{\Hom}{Hom}
\DeclareMathOperator{\Aut}{Aut}
\DeclareMathOperator{\rk}{rk}
\DeclareMathOperator{\ev}{ev}
\DeclareMathOperator{\Pic}{Pic}

\DeclareMathOperator{\supp}{supp}

\DeclareMathOperator{\Stab}{Stab}
\DeclareMathOperator{\Attr}{Attr}

\DeclareMathOperator{\codim}{codim}

\DeclareMathOperator{\Db}{D^b}

\newcommand{\vir}{\textup{vir}}

\newcommand{\pt}{\textup{pt}}

\newcommand{\QM}{\mathsf{QM}}

\newcommand{\cF}{\mathscr{F}}

\newcommand{\tO}{\widehat{\mathscr{O}}}

\newcommand{\Ami}{\cA_{\min}}

\DeclareMathOperator{\Coker}{Coker}
\DeclareMathOperator{\tr}{tr}

\DeclareMathOperator{\Spec}{Spec}

\DeclareMathOperator{\Coh}{Coh}
\DeclareMathOperator{\Casimir}{Casimir}

\newcommand{\St}{\mathbf{St}}

\newtheorem{Proposition}{Proposition}[section]
\newtheorem{Lemma}[Proposition]{Lemma}

\newtheorem{Theorem}{Theorem}

\theoremstyle{definition}
\newtheorem{Definition}{Definition}[section]

\begin{document}

\title{Enumerative geometry and geometric 
representation theory} 
\author{Andrei Okounkov}
\date{} 
\maketitle

\begin{abstract}
This is an introduction to: (1) the enumerative geometry of 
rational curves in equivariant symplectic resolutions, and 
(2) its relation to the structures of geometric representation 
theory. Written for the 2015 Algebraic Geometry
Summer Institute. 
\end{abstract}

%\setcounter{tocdepth}{2}
%\tableofcontents

\section{Introduction}

\subsection{}
These notes are written to accompany my lectures given in 
Salt Lake City in 2015.  
With modern technology, one should be able to access  
the materials from those lectures from anywhere in the world, 
so a text transcribing them is not really needed. Instead, I will 
try to spend more time on points that perhaps require too much
notation for a broadly aimed series of talks, but which should ease the 
transition to reading detailed lecture notes like \cite{PCMI}. 

The fields from the title are vast and they intersect in many different ways. 
Here we will talk about a particular meeting point, the progress at 
which in the time since Seattle I find exciting enough to report in 
Salt Lake City. The advances both in the subject matter itself, and 
in my personal understanding of it, owe a lot to M.~Aganagic, 
R.~Bezrukavnikov, P.~Etingof, D.~Maulik, N.~Nekrasov, and
others, as will be clear from the narrative.

\subsection{}
The basic question in representation theory is to 
describe the homomorphisms
\begin{equation}
\textup{some algebra $\cA$} \to 
\textup{matrices} \,,  \label{rtq} 
\end{equation}
and the geometric representation theory aims to
describe the source, the target, the map itself, or 
all of the above, geometrically. 

For example, matrices may be replaced by 
correspondences, by which we mean cycles in $X \times X$, 
where $X$ is an algebraic variety, or K-theory 
classes on $X \times X$ et cetera. These form 
algebras with respect to convolution and act 
on cycles in $X$ and $K(X)$, respectively. Indeed, 
$\Mat(n,R)$, where $R$ is a ring, is nothing but 
$K(X^{\times 2}) \otimes_\Z R$ where $X$ is a finite
set of cardinality $n$, and the usual rules of linear 
algebra, when written in terms of pullback, 
product, and pushforward, apply 
universally.  (This is also how e.g.\ integral operators
act.) 

Most of the time, these can and should be upgraded 
to Fourier-Mukai kernels \cite{Huy}, although doing this in 
the enumerative context should be weighted against the 
cost of losing deformation invariance --- a highly 
prized and constantly used property of 
pragmatically defined geometric counts. 

\subsection{} 

There is some base ring implicit in \eqref{rtq} and, in real life, 
this ring is usually $\Hd_G(\pt)$ or $K_G(\pt)$, where $G$ is a
reductive group acting on $X$ although, of course, generalizations 
are possible. 

In what follows, it will be natural and important to work with 
the maximal equivariance allowed by the problem.  Among 
their many advantages, equivariant counts are often defined when 
nonequivariant aren't. Concretely, the character of an
infinite-dimensional module $V$ may be well-defined as 
a rational function on $G$, but $1\in G$ will typically be a pole 
of this function.

\subsection{}\label{s_cM} 

Enumerative geometry is an endless source of 
interesting correspondences of the following kind. 
The easiest way for two distant points $x_1,x_2\in X$ 
to interact is to lie on a curve $C\subset X$ of some 
degree and genus and, perhaps, additionally constrained by 
e.g.\ incidence to a fixed cycle in $X$. 
Leaving the exact notion of a curve vague for a 
moment, one can contemplate a moduli space $\cM$ of 
two-pointed curves in $X$ with an evaluation map 
$$
\ev: \cM \to X \times X  \,.
$$
One can use this evaluation map to construct 
correspondences given, informally, by pairs of points 
$(x_1,x_2)\in X^{\times 2}$ 
that lie on such and such curve. For this, one needs
geometrically 
natural cycles or K-classes on $\cM$ to push forward
and, indeed, these are available in some 
generality. 

Deformation theory provides a local description of 
the moduli space $\cM$ near a given curve $C$ in 
$X$. Counting the parameters for deformations minus
the number of equations they have to satisfy, 
one computes 
\begin{equation}
\textup{expected} \, \dim \cM = 
(3 - \dim X) (g-1) + (d,c_1(X))+n \,,   \label{virdim} 
\end{equation}
where $g$ is the genus
of a curve $C\in \cM$, $d\in H_2(X,\Z)$ is its degree, 
and $n=2$ is the number of marked points. This is only a lower 
bound for the actual dimension of $\cM$ and it is seldom 
correct. Enumerative geometers treat this as they would treat 
any excess intersection problem: 
redundant equations still cut out a canonical 
cycle class $\left[\cM\right]_\vir$ of the correct dimension, known as the 
\emph{virtual fundamental cycle}, see \cite{BF}. 
There is a parallel
construction of the virtual structure sheaf 
$\cO_{\cM,\vir}\in K_{\Aut(X)}(\cM)$, see 
\cite{CF3,FGoe}.  As we deform $X$, the moduli spaces $\cM$ 
may jump wildly, but the curve 
counts constructed using $\left[\cM\right]_\vir$ and 
$\cO_{\cM,\vir}$ do not change. 

It will seem like a small detail now, but a certain symmetrized 
version 
\begin{equation}
\tO_\vir = \cO_\vir \otimes \cK_{\cM,\vir}^{1/2} \otimes \dots 
\label{tO1} 
\end{equation}
of the virtual structure sheaf has improved
self-duality properties and links better with both representation 
theory and mathematical physics. The square-root factor in 
\eqref{tO1} is the square root of the virtual canonical bundle, 
which exists in special circumstances \cite{NO}. The 
importance of such twist in enumerative K-theory was 
emphasized by Nekrasov \cite{Zth}, the rationale being that for K\"ahler 
manifolds the twist by a square root of the canonical bundle 
turns the Dolbeault operator into the Dirac operator. 
We will always use \eqref{tO1}, 
where the terms concealed by the dots will be specified after we 
specify the exact nature of $\cM$.

\subsection{}

One can talk about a representation-theoretic answer to 
an enumerative problem if the correspondences 
$\ev_* \left[\cM\right]_\vir$ or $\ev_* \tO_\vir$ are 
identified as elements of some sufficiently rich  
algebra acting by correspondences on $X$. Here rich 
may be defined pragmatically as allowing for computations or proofs. 

This is exactly the same as being able to place the 
evolution operator of a quantum-mechanical or 
field-theoretic problem into a rich algebra of operators acting 
on its Hilbert space. A mathematical physicist would call this 
phenomenon \emph{integrability}. Nekrasov and Shatashvili \cite{NS1,NS2} were
to first to suggest, in the equivalent language of supersymmetric 
gauge theories, that the enumerative problems discussed here are
integrable. 

\subsection{}

My personal intuition is that there are much fewer ``rich'' 
algebras than there are interesting 
algebraic varieties, which means there has 
to be something very special about $X$ to have a real link between 
curve-counting in $X$ and representation theory. In any case, 
the progress in enumerative geometry that will be described in these
lectures is restricted to certain very  special algebraic 
varieties. 

Ten years ago in Seattle, Kaledin already spoke about \emph{equivariant 
symplectic resolutions}, see \cite{Kal} and Section 
\ref{s_sympl_res}, and the importance of this 
class of algebraic varieties has been only growing 
since. In particular, it was understood by Bezrukavnikov
 and his collaborators 
that the geometry of rational curves in an equivariant symplectic 
resolution $X$ is tightly intertwined with the geometric 
structures described in \cite{Kal} and, specifically, 
with the derived autoequivalences of $X$ and the representation theory of 
its quantization. 

I find it remarkable that here one makes contact with a very 
\emph{different} interpretation of what it means for \eqref{rtq} to be 
geometric. A noncommutative algebra $\cX$ in the source of \eqref{rtq}
may be constructed geometrically as a \emph{quantization} of
a symplectic algebraic variety $X$ and this very much ties 
$\cX$-modules with coherent sheaves on $X$. We denote this 
noncommutative 
algebra $\cX$ partly to avoid confusion with some algebra $\cA$ 
acting on $X$ by correspondences but, in fact, in this subject, 
algebras of both kinds often trade places ! This means that 
$\cA$ turns out to be related to a quantization of some $X^\vee$, 
and vice versa. 

Dualities of this kind originate in supersymmetic 
gauge theories and are known under various names there, see e.g.\
\cite{IS,BHOO,BDG,BDGH}, mathematicians call them \emph{symplectic duality} 
following \cite{BLPW1}, see also \cite{BLPW2,NakC1,
NakC2,NakC3}. Perhaps at the next summer 
institute someone 
will present a definite treatment of these dualities. 

\subsection{}
The relation between rational curves in $X$  and 
representation theory of $\cX$ 
will be explained in due course
below; for now, one can note that perhaps a link between 
enumerative geometry and representation theory exists for a class
of varieties which is at least as large and about as special as 
equivariant symplectic resolutions. 

Prototypical examples of equivariant symplectic resolutions 
are the cotangent bundles $T^*G/P$ of projective homogeneous 
spaces\footnote{
Quantum cohomology of $T^*G/P$ was sorted out in \cite{BMO,Su}. 
Arguably, it is even simpler than the beautiful story we have 
for $G/P$ itself.}. Perhaps equivariant symplectic resolutions form the 
right generalization of semisimple Lie algebras for the needs
of today's mathematics ? 

% Certainly, one has to leave the world of symplectic resolutions 
% to find, for example, projective homogeneous varieties $G/P$. 
% These have a very beautiful enumerative geometry, as discovered
% by Givental and Kim and others \cite{}. The cotangent bundles
% $T^*G/P$ are the prototypical equivariant symplectic resolutions 
% and their zero sections $G/P$ may be described e.g.\  as
% the fixed points of a torus scaling the cotangent fibers. Laumon 
% spaces studied in \cite{} from the angle of enumerative geometry 
% have a similar description. 

While equivariant symplectic resolutions await
Cartan and Killing of the present day to classify them, the largest 
and richest class of equivariant symplectic resolutions known to date
is formed by Nakajima quiver varieties \cite{Nak94,Nak98}. These are associated 
to a quiver, which means a finite graph with possibly loops and multiple 
edges. Quivers generalize Dynkin diagrams (but the meaning 
of multiple edges is different).  {}From their very beginning, 
Nakajima varieties played a important role in geometric 
representation theory; that role has been only growing since. 

\subsection{}

Geometers particularly like Nakajima varieties associated to affine 
ADE quivers, because they are moduli of framed sheaves on the 
corresponding ADE surfaces. In particular, Hilbert schemes of 
points of ADE surfaces are Nakajima varieties. By definition, a map 
\begin{equation}
f: C \to \Hilb(S, \textup{points})\,,\label{fCtoHilb} 
\end{equation}
where $C$ is a curve and $S$ is a surface, is the same as a 
1-dimensional subscheme of $Y=C\times S$, flat over $C$, 
and similarly for other moduli of sheaves on $S$. 
This is the prosaic basic link between curves in affine ADE Nakajima 
varieties and the enumerative geometry of sheaves
on threefolds, known as the Donaldson-Thomas theory
\footnote{It is not unusual for two very different moduli 
spaces to have isomorphic open subsets, and it doesn't prevent 
the enumerative information collected from the two from 
having little or nothing in 
common. It is important for $S$ to be symplectic in order
for the above correspondence to remain uncorrected 
enumeratively in cohomology. A further degree of precision is 
required in K-theory.}. We are primarily interested in 1-dimensional 
sheaves on $Y$, and we will often refer to them as curves
\footnote{There is the following clear and significant difference between 
the curves in $\Hilb(S)$ and the corresponding curves in $Y$. The 
source $C$ in \eqref{fCtoHilb} is a \emph{fixed} curve, and from the 
enumerative point of view,  it may be taken to be a union of rational 
components. The ``curves'' in $Y$ being enumerated 
can be as arbitrary $1$-dimensional sheaves as our moduli spaces
allow.}. 

Other Nakajima varieties resemble moduli spaces of sheaves on 
a symplectic surface --- they can be interpreted as moduli of 
stable objects in certain 2-dimensional Calabi-Yau 
categories. Correspondingly, enumerative geometry of curves 
in Nakajima varieties has to do with counting stable objects
in 3-dimensional categories, and so belongs to the Donaldson-Thomas
theory in the broad sense. 

\subsection{}
To avoid a misunderstanding, a threefold $Y$ does not need to be 
Calabi-Yau to have an interesting enumerative geometry or 
enumerative K-theory of sheaves.  Calabi-Yau threefolds are certainly
distinguished from many points of view, but from the perspective of 
DT counts, either in cohomology or in K-theory, the geometry of 
curves in the good old projective space is at least as interesting. 
In fact, certain levels of complexity in the theory collapse under 
the Calabi-Yau assumption. 

When $Y$ is fibered over a curve in ADE surfaces, its DT theory 
theory may be directly linked to curves in the corresponding 
Nakajima variety. While this is certainly a very narrow class 
of threefolds, it captures the essential information about DT counts
in all threefolds in the following sense. 

Any threefold, together 
with curve counts in it, may be glued out of certain model 
pieces, and these model pieces are captured by the ADE fibrations 
(in fact, it is enough to take only $A_0$, $A_1$, and $A_2$
surfaces). This may be compared with Chern-Simons theory 
of real threefolds: that theory also has pieces, described by the 
representation theory of the corresponding loop group or 
quantum group at a root of unity, which have to be glued to obtain 
invariants \footnote{In fact, Donaldson-Thomas theory was 
originally developed as a complex analog of the Chern-Simons 
theory. It is interesting to note that while representation theory 
of quantum groups or 
affine Lie algebras is so-to-speak the DNA of Chern-Simons 
theory, the parallel role for Donaldson-Thomas theory is played
by algebras which are quantum group deformations of \emph{double}
loops.}. 

The algebraic version of breaking up a threefold $Y$ we need here
is the following. Assume that $Y$ may be degenerated, with a 
smooth total space of the deformation, to a transverse union 
$Y_1 \cup_D Y_2$ of two smooth threefolds along a smooth 
divisor $D$. Then, on the one hand, degeneration formula of 
Li and Wu \cite{LiWu} gives
\footnote{In K-theory, there is a correction to the gluing 
formula as in the work of Givental \cite{Giv1}.}
 the curve counts in $Y$ in terms of counts
in $Y_1$ and $Y_2$ that also record the \emph{relative} 
information at the divisor $D$. On the other hand, Levine and 
Pandharipande prove \cite{LP} that the relation 
\begin{equation}
\left[Y\right] = \left[Y_1\right]  +  \left[Y_2\right] - 
\left[\bP(\cO_D\oplus N_{Y_1/D})\right] \label{algcob} 
\end{equation}
generates all relations of the algebraic cobordism, and so, in
particular, any projective $Y$ may be linked to a product of projective 
spaces by a sequence of such degenerations
\footnote{This does not trivialize DT counts, in particular, does not 
make them factor through algebraic cobordism, for the following
reasons. 

First, degeneration formula glues \emph{nonrelative}
(also known as absolute) DT counts out of \emph{relative}
 counts in $Y_1$ and $Y_2$ and to use it again one needs 
a strategy for replacing relative conditions by absolute ones. 
Such strategy has been perfected by Pandharipande and Pixton
\cite{PP1,PP2}. It involves certain universal substitution rules, which can
be, again, studied in model geometries and are best understood
in the language of geometric representation theory \cite{Sm1}. 

Second, one doesn't always get to read the equation \eqref{algcob} 
from left to right: sometimes, one needs to add a component to 
proceed. This means one will have to solve for curve counts in $Y_1$ 
in terms of those in $Y$ and $Y_2$, and a good strategy for this 
is yet to be developed.} . 

If $Y$ is toric, e.g.\ a product of projective spaces, then 
equivariant localization may be used to break up curve counts 
into further pieces and those are all captured by 
$Y=\bP^1 \times A_n$ geometries, see \cite{MOOP}.

\subsection{}\label{sHopf} 

Nakajima varieties $X$ are special even among equivariant symplectic 
resolutions, in the following sense. The algebras $\cA$ acting 
by correspondences on $X$ are, in fact, \emph{Hopf algebras}, which 
means that there is a coproduct, that is, an algebra homomorphism 
$$
\Delta: \cA \to \cA \otimes \cA \,,
$$
a counit $\varepsilon: \cA \to \textup{base ring}$, 
 and an antipode, satisfying standard axioms (a certain completion 
may be needed for infinite-dimensional $\cA$). More abstractly, these
mean that the category $\cC=\cA\textup{-mod}$ has 
\begin{itemize}
\item a tensor product $\otimes$, which need \emph{not} be commutative, 
\item a trivial representation, given by $\varepsilon$, which is identity
  for $\otimes$, 
\item left and right duals, compatible 
  with $\otimes$. 
\end{itemize}
The notion of a tensor category is, of course, very familiar to 
algebraic geometers, except possibly for the part that allows for 
\begin{equation}
V \otimes W \not\cong W \otimes V \,.
\label{nctens} 
\end{equation}
This noncommutativity is 
what separates representation theory of \emph{quantum groups} from 
 representation theory of 
usual groups. Representation-theorists know that a mild 
noncommutativity makes a tensor category even more 
constrained and allows for an easier reconstruction of $\cA$ from
$\cC$, see e.g.\ \cite{Ebook1,Ebook2} for an excellent exposition. 

The antipode, if it exists, is uniquely reconstructed from the rest of 
data (just like the inverse in a group). We won't spend time on it 
in these notes.

\subsection{}

Specifically, in K-theoretic setting, we will have
$\cA=\cU_\hbar(\gh)$ for a certain Lie algebra $\fg$. 
This means that 
$\cA$ is a Hopf algebra deformation\footnote{with 
deformation parameter $\hbar$, which is the equivariant 
weight of the symplectic form on $X$ in its geometric origin} 
of $\cU(\fg [t^{\pm1}])$, where $\fg [t^{\pm1}]$ denotes the Lie 
algebra of Laurent polynomials with values in $\fg$. 

The Lie algebra $\fg$ is 
by itself typically infinite-dimensional and even, strictly 
speaking, 
infinitely-generated. The action of $\cA$ will extend 
the action of $\cU_\hbar(\gh_{\textup{\sc km}})$ for a certain 
Kac-Moody Lie algebra $\fg_{\textup{\sc km}} \subset \fg$ 
defined by Nakajima via an explicit assignment on generators
\cite{Nak98}. 

The noncommutativity \eqref{nctens} will be mild and, in its 
small concentration, very beneficial for the development of the theory. 
Concretely, the ``loop rotation'' automorphism 
\begin{equation}
t \mapsto u\cdot t \,, \quad u \in \Ct \,, \label{loop_rot} 
\end{equation}
of $\cU(\fg [t^{\pm1}])$ deforms, and so every $\cA$-module 
$V$ gives rise to a family of modules $V(u)$ obtained by 
precomposing with this automorphism.  As it will turn out, 
the commutativity in \eqref{nctens} is restored after a 
generic shift, that is, there exists an $\cA$-intertwiner 
$$
R^\vee(u_1/u_2): V_1(u_1) \otimes V_2(u_2) 
\xrightarrow{\,\,\sim\,\,} 
V_2(u_2) \otimes V_1(u_1)
$$
which is a rational 
function of $u=u_1/u_2\in \Ct$. The poles and zeros of 
$\det R^\vee$ correspond to the values of the parameter for which 
the two tensor products are really not isomorphic. 

\subsection{}\label{s_introR}

Usually, one works with the 
so-called $R$-matrix
$$
R(u) = (12) \cdot R^\vee \in \End (V_1 \otimes V_2) \otimes 
\Q(u) \,, 
$$
which intertwines the action defined via the 
coproduct $\Delta$ and via the opposite coproduct 
$\Delta_\textup{opp}=(12)\Delta$. We denote by
$\otimes_\textup{opp}$ this opposite way to take 
tensor products. 

The shortest way to reconstruction of $\cA$ from $\cC$ is 
through the matrix coefficients of the $R$-matrix. 

Concretely, since $R$ is an operator in a tensor product, its matrix 
coefficients in $V_1$ give  $\End(V_2)$-valued functions of $u$. 
The coefficients of the expansion of these functions 
around $u=0$ or $u=\infty$ give us countably many operators 
in $V_2$, which starts looking like an algebra of the right size to be 
a deformation of $\cU(\fg [t^{\pm1}])$. The Yang-Baxter 
equation satisfied by $R$ 
\begin{multline}
R_{V_1,V_2 }(u_1/u_2) \, R_{V_1,V_3}(u_1/u_3) \, R_{V_2,V_3}(u_2/u_3) =\\=
R_{V_2,V_3}(u_2/u_3) \, R_{V_1,V_3}(u_1/u_3) \, R_{V_1,V_2}(u_1/u_2)\,, \label{YB} 
\end{multline}
which is an equality of rational functions with values in 
$\End(V_1 \otimes V_2 \otimes V_3)$, is a commutation 
relation for these matrix elements. 

Thus, a geometric construction of an $R$-matrix gives 
a geometric construction of an action of a quantum 
group.

\subsection{}
Nakajima varieties, the definition of which will be recalled below, 
are indexed by a quiver and two dimension vectors 
$\bv,\bw \in \N^{\textup{vertices}}$ and one sets 
\begin{equation}
X(\bw) = \bigsqcup_\bv X(\bv,\bw)\,. \label{Xw} 
\end{equation}
For example, for framed sheaves on ADE surfaces, $\bw$ is the data
at infinity of the surface 
(it includes the rank), while $\bv$ records possible values of
Chern classes. In the eventual description of $K_\bT(X(\bw))$ as an 
module over a quantum group, the decomposition \eqref{Xw} will 
be the decomposition into the weight spaces for a Cartan 
subalgebra $\fh \subset \fg$. 

To make the collection $\{K_\bT(X(\bw))\}$ into modules over 
a Hopf algebra, one needs tensor product. The seed from which 
it will sprout is very simple and formed by an inclusion 
\begin{equation}
X(\bw) \times X(\bw') \hookrightarrow X(\bw+\bw') 
\label{Xwww} 
\end{equation}
as a fixed locus of a torus $\bA \cong \Ct$. For framed sheaves, 
this is the locus of direct sums. The pull-back in K-theory under
\eqref{Xwww}, while certainly a worthwhile map, is not what we 
are looking for because it treats the factors 
essentially symmetrically. 

Instead, we will use a certain canonical correspondence 
$$
\Stab_{\fC,\dots}: K_\bT(X^\bA)  \to K_\bT(X)\,, 
$$
called the \emph{stable envelope}. 
It may be defined in some generality for a pair of tori 
$\bT\supset\bA$ acting on an algebraic
 symplectic manifold $X$ so that 
$\bA$ preserves the symplectic form. As additional 
input, the correspondence takes a certain cone $\fC \subset \Lie\bA$,
and other ingredients which will be discussed below. 
For $\bA\cong\Ct$, the choice of $\fC$ is a choice of $\R_{\gtrless 0} 
\subset \R$ or, equivalently, the choice of $\otimes$ versus 
$\otimes_\textup{opp}$. 

We demand that the two solid arrows in the diagram 
\begin{equation}
\xymatrix{
K_\bT(X(\bw))(a) \otimes K_\bT(X(\bw')) \ar[rrd]^{\Stab_+}
\ar@{-->}[dd]_{R(a)}\\
&& K_\bT(X(\bw+\bw')) \\
K_\bT(X(\bw))(a) \otimes_\textup{opp} K_\bT(X(\bw')) \ar[rru]_{\Stab_-}
}\label{SSR} 
\end{equation}
are morphisms in our category. By $K_\bT(\pt)$-linearity, 
this forces the rational 
vertical map $R(a)$ because stable envelopes are isomorphisms 
after localization. 

Here $a\in \bA \cong \Ct$ is an automorphism 
which does not act on $X^\bA$ but does act on vector bundles on $X^\bA$
through their linearization inherited from $X$.  This twists the 
module structure by an automorphism as in Section \ref{sHopf}
and identifies $a$ with the parameter $u=u_1/u_2$ in the R-matrix. 

As already discussed, the assignment of $R$-matrices 
completely reconstructs the Hopf algebra with its representations. 
The consistency of this procedure, in particular the Yang-Baxter 
equation \eqref{YB}, follow from basic properties of stable envelopes. 

\subsection{}
One may note here that even if one is primarily interested in 
Hilbert schemes of ADE surfaces, it is very beneficial to think 
about framed sheaves of higher rank. In higher rank, we have nontrivial 
tori $\bA$ acting by changing of framing and hence the associated 
$R$-matrices.  

Similarly, it will prove easier to go through enumerative 
K-theory of curves in the moduli spaces of higher-rank sheaves
even if one is only interested in curves in the Hilbert scheme of
points. 

At this point, we have introduces the basic curve-counting 
correspondences in $X$ and certain algebras $\cA$ 
acting by correspondences in $X$. The main goal of these note
is to discuss how one identifies the former inside the latter. 

\subsection{Acknowledgments}

\subsubsection{}
This is a report about a rapidly developing field. This is 
certainly exciting, with the flip side that the field itself, 
let alone its image in the author's head, is very far 
from equilibrium. What I understand about the field it is very 
much a product of my interactions with M. Aganagic, R. Bezrukavnikov,
P. Etingof, D. Maulik, I.~Losev, 
H.~Nakajima, N. Nekrasov, and others, both within 
the format of a joint work and outside of it. I am very 
grateful to all these people for what I have learned 
in the process. 

\subsubsection{}
I am very grateful to the organizers of the Salt Lake City institute
for creating this wonderful event and for 
the honor to address the best algebraic geometers of the 
universe. Many thanks to the AMS, the Simons Foundation, 
the Clay Mathematics Institute, the NSF and other funding 
agencies that made it possible. 

I am personally deeply 
grateful to the Simons Foundation for being supported in the 
framework of the Simons Investigator program. 
Funding by the Russian Academic Excellence Project '5-100' is 
gratefully acknowledged.

\section{Basic concepts}

\subsection{Symplectic resolutions}\label{s_sympl_res}

\subsubsection{}

By definition \cite{Beau,Kal}, a \emph{symplectic resolution} is a smooth 
algebraic symplectic variety $(X,\omega)$ such that the 
canonical map 
$$
X \xrightarrow{\quad}  X_0 = \Spec H^0(\cO_X) 
$$
is projective and birational. 

For us, the focus is  the smooth variety $X$ and the singular 
affine variety $X_0$ plays an auxiliary role, because, for example, 
there are no curves to count in $X_0$. Of course, there exists also an 
exactly complementary point of view that singularities are 
essential and resolutions --- auxiliary. 

\subsubsection{} 

A symplectic resolution $X$ is called \emph{equivariant} or conical 
if there is a $\Ct$ action on $X$ that scales $\omega$ and 
contracts $X_0$ to a 
a point. In other words, there is a grading on $\C[X_0]$ such that 
$$
\C[X_0]_d =
\begin{cases}
\C\,, & d=0 \,,\\
0 \,, & d<0 \,. 
\end{cases}
$$
The combination of the two requirements implies $\Ct$ must scale
$\omega$ by a nontrivial character and hence 
$$
[\omega] = 0 \in H^2(X,\C) \,, 
$$
since $\Hd(X)$ is a trivial $\Aut(X)_\textup{connected}$-module. 

\subsubsection{} 

The prototypical example of an equivariant symplectic resolution 
is $T^*G/P$ where $P\subset G$ is a parabolic subgroup of a 
semisimple Lie group. In particular, 
\begin{equation}
T^*G/B \xrightarrow {\quad \textup{moment map}\quad } 
\textup{Nilcone} \subset \mathfrak{g}^* \cong \mathfrak{g}
\label{Nilcone} 
\end{equation}
is a famous resolution with many uses in geometric 
representation theory. 

For an example with a less pronounced Lie-theoretic flavor one 
can take 
$$
\Hilb(\C^2, \textup{$n$ points}) 
\xrightarrow{\quad \textup{Hilbert-Chow} \quad}  S^n \C^2 \,,
$$
where the symplectic form comes from the standard 
$\omega_{\C^2} = dx_1 \wedge dx_2$. This is a basic example 
of a Nakajima variety. 

\subsection{Nakajima quiver varieties}

\subsubsection{} \label{s_Nak} 
Nakajima quiver varieties \cite{Nak94,Nak98} are defined as algebraic 
symplectic reductions 
\begin{align}
X &= T^* M \rd G_\bv \notag \\
 & = \textup{moment map}^{-1}(0) \rdd_\theta G_\bv \,,
\label{defNak} 
\end{align}
where $G_\bv=\prod GL(V_i)$ and $M$ is a linear representation of 
a special kind. Namely, $M$ is 
a sum of the defining and so-called bifundamental 
representations, that is, 
\begin{equation}
M = \bigoplus_i \Hom(W_i,V_i) \oplus \bigoplus_{i\le j} 
\Hom(V_i,V_j) \otimes Q_{ij} \,, \label{descrM} 
\end{equation}
where $W_i$ and $Q_{ij}$ are certain multiplicity spaces on which 
$G_\bv$ does not act. This gives an action 
\begin{equation}
\prod GL(W_i) \times \prod GL(Q_{ij}) \times \Ct_\hbar \to 
\Aut(X) \,, \label{GLGL} 
\end{equation}
where the last factor scales the cotangent directions with 
weight $\hbar^{-1}$. This inverse is conventional, it gives weight 
$\hbar$ to the symplectic form on $X$. 

\subsubsection{} 
In the GIT quotient in \eqref{defNak}, 
the stability parameter $\theta$ is a character of $G_\bv$, up 
to proportionality. It has to avoid the walls of a certain finite central 
hyperplane arrangement for the quotient to be smooth. 

What is special about representations \eqref{descrM} is that 
the stabilizer in $G_\bv$ of any point in $T^*M$ is cut out 
by linear equations on matrices and is the set of 
invertible elements in an associative algebra over $\C$. 
Therefore, it cannot be a nontrivial finite group. 
In general, algebraic symplectic reductions 
output many Poisson orbifolds, but it is difficult to 
produce a symplectic resolution in this way (or in 
any other way, for that matter). 

\subsubsection{}
The data of a group $G_\bv$ and a
 representation \eqref{descrM} is encoded by 
a quiver $Q$ and two dimension vectors. One 
sets
\begin{align*}
I &= \textup{vertices of $Q$} \\
 &= \textup{factors in $G_\bv$} 
\end{align*}
and joins two vertices $i,j\in I$ by $\dim Q_{ij}$ edges. The vectors
$$
\bv = (\dots, \dim V_i, \dots)\,, \quad 
\bw= (\dots, \dim W_i, \dots) 
$$
belong to $\N^I$ which makes them dimension vectors in 
quiver terminology. 

Nakajima varieties $X(\bv,\bw)$ corresponding to the same 
quiver $Q$ have a great deal in common and it is natural to group them 
together. In particular, it is natural to define $X(\bw)$ as in 
\eqref{Xw}. 

\subsubsection{}
The apparent simplicity of the definition \eqref{defNak} is
misleading. For example, let $Q$ be the quiver with one vertex and 
one loop. Then $X(\bv,\bw)$ is the moduli space of torsion-free
sheaves $\cF$ on $\bP^2$, framed along a line, with 
$$
(\rk \cF, c_2(\cF) )  = (\bw,\bv) \in \N^2 \,. 
$$
This variety has deep geometry and plays a very 
important role in mathematical physics, in particular, in the study 
of 4-dimensional supersymmetric gauge theories and instantons.  

Nakajima varieties have other uses in supersymmetric gauge 
theories as they all may be interpreted as certain moduli of 
\emph{vacua}, which is distinct from their interpretation as 
\emph{instanton} moduli for quivers of affine ADE type. 

\subsubsection{}\label{stensNak}

Let a rank $n$ torus $\bA$ with coordinates $a_i\in\Ct$ act
on $X(\bw)$ via 
$$
W_k = \bigoplus_{i=1}^n a_i \otimes W_k^{(i)} \,,
$$
where $W_k^{(i)}$ are trivial $\bA$-modules.  The fixed points of 
this action correspond to analogous gradings 
$$
V_k = \bigoplus_{i=1}^n a_i \otimes V_{k}^{(i)} \,,
$$
and to arrows that preserve it. Therefore 
$$
X(\bw)^\bA = \prod_{i=1}^n X(\bw^{(i)}) \,, 
$$
with the inherited stability condition. 
The case $n=2$ with $(a_1,a_2) = (a,1)$ appears in \eqref{Xwww}.

\subsection{Basic facts about rational curves in $X$}

\subsubsection{}

Following the outline of section \ref{s_cM}, but still not sewing 
permanently was tacked there, consider a moduli 
space $\cM$ of rational curves with 2 marked points 
in an equivariant symplectic resolution $X$. {}From formula
\eqref{virdim}, the virtual dimension of $\cM$ equals 
$$
\vir \dim \cM = \dim X - 1 \,,
$$
independent of the degree, because $c_1(X)=0$.
 This will be the dimension of 
$\ev_* \left[\cM\right]_\vir$ in $X \times X$. 

\subsubsection{}\label{s_deform} 

Deformations  of $(X,\omega)$ are described\footnote{
For Nakajima varieties, these deformations may be described 
explicitly by changing the value of moment map in \eqref{defNak}.}
by the period map 
$$
[\omega] \in H^2(X,\C) \,, 
$$
see \cite{Kal}. In particular, a generic deformation $X'$ will have 
no classes $\alpha \in H_2(X',\Z)$ such that $\int_\alpha \omega_{X'} =0$ 
and hence no algebraic curves of any kind (in fact, $X'$ will be
affine). Therefore any deformation-invariant curve counts in $X$ 
must vanish. 

To avoid such trivialities, we consider equivariant curve counts, 
where we must include the action of the $\Ct_\hbar$ factor that 
scales the symplectic form. It scales nontrivially the
identification of $H^2(X,\C)$ with the base of the deformation, 
hence there are no $\hbar$-equivariant deformations of $X$. 

Reflecting this, we have
$$
\left[\cM\right]_\vir = \hbar \,
\left[\cM\right]_{\vir,\textup{reduced}}
$$
where $\hbar \in \Hd_{\Aut(X)}(\pt)$ is the weight of the symplectic 
form and 
$$
\dim \left[\cM\right]_{\vir,\textup{reduced}} = \dim X \,. 
$$

\subsubsection{}
Clearly, 
$$
\ev(\cM) \subset X \times_{X_0} X 
$$
and by a fundamental property of symplectic resolutions, see
\cite{Kaledin_Poisson, Nami}, 
\begin{equation}
\St= X \times_{X_0} X \subset X \times X
\label{def_Stein} 
\end{equation}
is at most half-dimensional with Lagrangian half-dimensional
components. In the example \eqref{Nilcone}, the subscheme 
\eqref{def_Stein} is known as the \emph{Steinberg} 
variety, and it is convenient to extend this usage to general 
$X$. 

Putting two and two together, we have the 
following basic 

\begin{Lemma}[\cite{BMO}] \label{l1} 
For any degree, $\ev_* \left[\cM\right]_{\vir,\textup{reduced}}$ is a
$\Q$-linear combination of Lagrangian components of the Steinberg 
variety. 
\end{Lemma}

Rational coefficients appear here because
$\left[\cM\right]_{\vir}$ is defined with 
$\Q$-coefficients, due to automorphisms of objects 
parametrized by $\cM$.  Intuitively, it seems unlikely that that 
the coefficients in Lemma \eqref{l1} really depend on the 
details of the construction of $\cM$, as long as it has a perfect
obstruction theory. In fact, for Nakajima varieties, all possible 
flavors of curve counting theories give the same answer in 
cohomology. 

\subsubsection{}
One can make a bit more progress on basic principles. 
If $d\in H_2(X,\Z)_{\textup{eff}}$ is a degree of a rational 
curve, we denote by $z^d$ the corresponding element of the 
semigroup algebra of $H_2(X,\Z)_{\textup{eff}}$. The spectrum 
of this semigroup algebra is an affine toric variety; it is 
a toric chart corresponding to $X$ is the so-called 
K\"ahler moduli space $\fK$, see e.g.\ \cite{CoxKatz} for a discussion of the 
basic terminology of quantum cohomology. 

For a Nakajima variety $X$, 
$$
\fK =  \overline{H^2(X,\C)/H^2(X,\Z)} 
$$
 is a toric variety corresponding to the 
fan  in $H^2(X,\Z)$ formed by the ample cones of all flops of $X$. 
By surjectivity theorem of \cite{mn},
these are the cones of nonsingular 
values of the stability parameter $\theta$. 

We consider the operator 
of \emph{quantum multiplication} in $\Hd_\bT(X)$ by 
$\lambda\in H^2(X)$ 
\begin{equation}
\lambda \, \star \, \cdot \, = \lambda \cup \, \cdot \, + \hbar \, 
\sum_{d} (\lambda, d)\,  z^d \, \ev_*
\left[\cM_d\right]_{\vir,\textup{reduced}} \,,  \label{qmult} 
\end{equation}
where $\cM_d$ is the moduli space of curves of degree $d$. 
Here $\lambda \cup \, \cdot \,$ is the operator of cup 
product by $\lambda$, it is supported on the diagonal of 
$X \times X$ as a correspondence. On very general grounds, 
the operators \eqref{qmult} commute and, moreover, the 
quantum, or Dubrovin, connection with operators 
\begin{equation}
\nabla_\lambda = \varepsilon \frac{\partial}{\partial \lambda} - 
\lambda \, \star \, \cdot \, \,, \quad 
\frac{\partial}{\partial \lambda} z^d = (\lambda,d) \, z^d
\label{qconn} 
\end{equation}
is flat for any $\varepsilon$. 

For equivariant symplectic resolutions, one conjectures that
after a shift 
\begin{equation}
z^d \mapsto (-1)^{(\bk_X,d)} z^d \label{kappa_shift} 
\end{equation}
by a certain canonical element $\bk_X\in H^2(X,\Z/2)$
called the canonical theta characteristic in \cite{MO}, the 
quantum connection takes the following form: 
\begin{equation} 
\textup{quantum part in \eqref{qmult}} 
\overset ? =
\hbar 
\sum_\alpha (\lambda,\alpha)\, \frac{z^\alpha}{1-z^\alpha}  \,
L_\alpha + \dots \,, \quad L_\alpha \in H_\textup{top}(\St) \,, 
\label{polesqmult} 
\end{equation}
where dots stand for a multiple of the diagonal component of $\St$,
that is, for a scalar operator.

Such scalar ambiguity is always
present in this subject and is resolved here by the fact that 
$$
\forall \gamma\in \Hd(X)\,, \quad 1 \star \gamma = \gamma\,, 
$$
and so the operator \eqref{polesqmult} annihilates $1\in H^0(X)$. 
The sum \eqref{polesqmult} is over a certain finite set of 
effective classes $\alpha$. 

It is clear from \eqref{polesqmult} that 
\eqref{qconn} is a connection on $\fK$ with 
regular singularities.  
Conjecture \eqref{polesqmult} is known for all concrete $X$ discussed
in these 
notes\footnote{Concretely, it is known for Nakajima varieties by 
\cite{MO}, for $T^*G/P$ by \cite{Su}, and for hypertoric 
varieties by \cite{McBSh}.} and will be assumed in what follows. 

Let $\Ami$ be the algebra of endomorphisms of 
$\Hd_\bT(X)$ generated by the operators of cup product and 
$H_\textup{top}(\St)$. Clearly, the operators 
\eqref{qmult} lie in $\Ami[[z]]$. In all analyzed examples, these have a
simple spectrum for generic $z$ and thus generate the full 
algebra of operators of quantum product 
by $\gamma \in \Hd_\bT(X)$; one expects this to happen in 
general. It would mean \eqref{qmult} generate a family, 
parametrized by $z$, of maximal commutative subalgebras of 
$\Ami$ which deforms the algebra of cup products for $z=0$. 

In examples, this matches very nicely with known ways to 
produce maximal commutative subalgebras. For instance, 
for Nakajima varieties, the operator \eqref{qmult} lie in 
what is known as \emph{Baxter} maximal 
commutative subalgebra in a certain $\cA$ as a corollary of the 
main result of \cite{MO}. 

One can imagine the above 
structures to be extremely constraining for general 
equivariant symplectic resolutions. 

\subsubsection{}

Now, in K-theory, the only thing that carries over from the 
preceding discussion is that, for any group $G$, 
\begin{equation}
\ev_* : 
K_G(\cM) \to 
K_G (\St) \,,  \label{evK} 
\end{equation}
where there right-hand side denotes K-classes with support in 
$\St\subset X^{\times 2}$. This is an algebra under convolution. 

The best case scenario is when one has a complete control 
over $K_G (\St)$. For example, one of the main results of 
\cite{CG}, see Chapter 7 there, shows 
$$
K_{G\times \Ct_\hbar}(\St) = \textup{affine Hecke algebra}\,,
$$
for $X=T^*G/B $. This is 
a beautiful algebra with a beautiful presentation; perhaps it would 
be too much to expect an equally nice description in general. 
The approach explained in these notes doesn't assume 
any knowledge about $K_G(\St)$ as a prerequisite. 

\subsubsection{}

Even when the target of the map \eqref{evK} is under control, 
it is not so easy to identify, say, 
$\ev_* \cO_{\cM,\vir}$ in it. The dimension argument, which 
identifies $\ev_* \left[ \cM \right]_{\vir,\textup{reduced}}$ with
linear combinations of component of $\St$ obviously does not apply in 
K-theory, so new ideas are needed. 
 There exists a body of work, notably by Givental 
and his collaborators \cite{Giv2,Giv3,Giv4}, that aims to lift the 
general structures of quantum cohomology to quantum 
K-theory of an equally general algebraic variety $X$. This 
theory has been studied for $G/B$ itself and also, for example, 
for toric varieties.  It would be interesting to see what 
it outputs for $T^*G/B$. 

A different strategy  for dealing with the complexities of 
K-theory works only for special $X$ and 
is based on a systematic use of rigidity and self-duality 
arguments, see for example \cite{PCMI} for an introduction.  For them 
to be available, one needs the obstruction theory of $\cM$ to 
have a certain degree of self-duality in the first place and to work 
with symmetrized K-classes on $\cM$ such as \eqref{tO1}. 

\subsubsection{}
Let 
$$
f: (C,p_1,p_2) \to X
$$
be a map of a 2-pointed rational curve to $X$, which we assume 
is symplectic, and so its tangent bundle $TX$ is self-dual up-to 
the weight of the symplectic form 
$$
(TX)^\vee = \hbar \, TX \,. 
$$
This makes the obstruction theory 
$$
\Def f - \Obs f = \Hd(C, f^* TX) 
$$
self-dual up to a certain correction 
\begin{equation}
(\Def-\Obs)^\vee +\hbar \, (\Def-\Obs) = 
\hbar\, \Hd(C,f^*TX\otimes(\cO_C - \cK_C)) \,. 
\label{dual_corr} 
\end{equation}
For an irreducible curve, we have  $\cO_C - \cK_C = \cO_{p_1} + \cO_{p_2}$, 
equivariantly with respect to $\Aut(C,p_1,p_2)$, and we can use 
marked points to compensate
\footnote{
These contributions at the marked points are the dots in 
\eqref{tO1}, see \cite{PCMI} for more details.}
 for the correction in \eqref{dual_corr}. 
But this breaks down when $C$ is allowed to break and we get 
$$
\cO_C - \cK_C \ne \cO_{p_1} + \cO_{p_2}
$$
away from the locus of chains of rational curves joining $p_1$ and 
$p_2$. 

The moduli space $\overline{\cM}_{0,2}(X)$ of stable rational maps to
$X$ allows arbitrary trees as domains of the map, and so it 
appears problematic to make its virtual structure sheaf 
self-dual. As an alternative, one can use the moduli space 
$\QM(X)$ of stable quasimaps to $X$. This has its origin in 
supersymmetric gauge theories and 
and is known to geometers in various flavors. The version 
best suited for our needs has been designed by Ciocan-Fontanine, 
Kim, and Maulik in \cite{CKM}. It requires a GIT presentation 
\begin{equation}
X = \widetilde{X} \rdd G \label{XGIT} 
\end{equation}
where $G$ is reductive and $\widetilde{X}$ is an affine algebraic
variety
with at worst locally complete 
intersection singularities in the unstable locus. It remains to be 
seen which of the equivariant symplectic resolutions can be 
presented in this way; for Nakajima varieties, it is provided by  
their construction. When available, quasimaps spaces have 
marked advantages in enumerative K-theory. 

\subsubsection{}

For example, 
for $X=\Hilb(\C^2,\textup{points})$, the quasimaps from a fixed 
curve $C$
to $X$ are naturally identified with stable pairs or points in 
Pandharipande-Thomas moduli spaces \cite{PT} for 
$Y=C\times \C^2$. A stable pair is a 
pure 1-dimensional sheaf $\cF$ on $Y$ with a section 
$$
\cO_Y \xrightarrow{\,\, s \,\,} \cF 
$$
such that $\dim \Coker s = 0$. These moduli spaces have 
several advantages over other DT moduli spaces, including
the Hilbert schemes of curves, and are used very 
frequently.  Marked points $p_1,p_2\in C$ introduce 
a \emph{relative} divisor 
$$
D = D_1 \sqcup D_2 \subset Y 
$$
where $D_i$ is the fiber over $p_i$. To have an evaluation map 
at these marked points, the curve $C$ is allowed to develop a 
chain of rational components at each $p_i$; note that the whole 
curve $C$ remains a chain in this process. To have a moduli space 
of maps from a nonrigidified rational curve, we must quotient 
the above by $\Aut(C,p_1,p_2)$. These are known affectionately as 
\emph{accordions} among the practitioners. 

Quasimaps to a general Nakajima variety behave in a very 
similar way, see \cite{PCMI} for an elementary introduction. 

\subsubsection{}

In fact, for quasimaps from either a rigid 2-pointed curve $C$ or 
from pure accordions to a Nakajima variety $X$, one gets
 the \emph{same} answer 
 \begin{equation}
\bfG = \ev_* \tO_{\QM(X)}  \in K_{\Aut X}(\St) \label{defglue} 
\end{equation}
known as the \emph{glue} matrix, see Theorem 7.1.4 in \cite{PCMI}. 
There is also a $q$-difference connection generalizing 
\eqref{qconn}. The geometric meaning of the 
increment $q$ in this $q$-difference 
equation is 
$$
q \in \Ct = \Aut(C,p_1,p_2) \,, \quad (C,p_1,p_2) \cong
(\bP^1,0,\infty)\,. 
$$
In fact, the variable $\varepsilon$ in \eqref{qconn} belongs 
to $\Lie \Aut(C,p_1,p_2)$ from its geometric origin. 
In contrast to cohomology, the operator of the 
quantum $q$-difference equation is much more involved than 
\eqref{defglue}. In particular, the matrix of the $q$-difference 
equation depends on $q$, whereas the glue matrix 
doesn't. The glue matrix \eqref{defglue} may be 
obtained from the quantum difference equation by a certain 
limit. 

Quantum difference equations for all Nakajima varieties have
been determined in \cite{OS}. The description is
representation-theoretic 
and a certain language needs to be developed before we can 
state it.  

\section{Roots and braids}

\subsection{K\"ahler and equivariant roots} 

\subsubsection{}

If equivariant symplectic resolutions are indeed destined to generalize 
semisimple Lie algebras, one should be able to say what 
becomes of the classical root data in this more general 
setting. 

For a semisimple Lie group  $G$, we have: 
\begin{itemize}
\item[---] roots, which are the nonzero weights of a 
maximal torus $\bA\subset G$ in the adjoint representation, 
\item[---] coroots, which have to do with special maps 
$SL(2) \to G$. 
\end{itemize}

The resolution $X=T^*G/B$ by itself does not distinguish
between locally isomorphic groups.  It is the adjoint group 
$G_\textup{ad}$ that acts naturally by symplectic automorphisms
of $X$. On the other hand, 
$$
H_2(X,\Z) = \textup{coroot lattice of $G_\textup{sc}$} \,, 
$$
where $G_\textup{sc}$ is the simply-connected group 
and this lattice is spanned by the images
$$
\bP^1 = SL(2)/{B_{SL(2)}}  \xrightarrow{\,\,  
\textup{coroot} \,\, } G/B \subset X \,. 
$$
The equivariant and K\"ahler roots of symplectic resolutions will generalize
the roots and coroots for $X=T^*G/B$, respectively. 
 They capture the weights 
of torus action and special rational curves in $X$, respectively. 

% They will naturally be of adjoint and simply-connected flavor, 
% respectively, on the two sides. 

\subsubsection{}

Let 
$$
\bA \subset \Aut(X,\omega) 
$$
be a maximal torus in the group of symplectic 
automorphisms of an equivariant symplectic resolution $X$. 

\begin{Definition}
The \emph{equivariant} roots of $X$ are the weights of a maximal 
torus $\bA $ in the normal bundle 
$N_{X/X^\bA}$ to its fixed locus $X^\bA$. 
\end{Definition}

For $X=T^*G/B$, this gives the roots of the adjoint group. 
For $X=\Hilb(\C^2,n)$, the torus $\bA$ is the maximal torus 
in $SL(\C^2)$.  We have 
$$
X^\bA = \{\textup{monomial ideals} \} 
$$
and the weights of $\bA$ in the tangent space to these points
are classically computed in terms of the hook-lengths for the 
corresponding partitions of $n$. Therefore 
$$
\textup{equivariant roots of $\Hilb(\C^2,n)$} = \{\pm1, 
\dots,\pm n\} \,.
$$
We see that, in contrast to finite-dimensional Lie theory, roots may be 
proportional for symplectic resolutions.

\subsubsection{}\label{scones}

As in the classical Lie theory, equivariant roots define root
hyperplanes in $\Lie \bA$ and these partition the real locus 
into a finite set of chambers $\fC$.  Once a chamber is 
fixed, it divides equivariant roots into positive and 
negative. 

Geometrically, this means splitting normal directions to 
$X^\bA$ into attracting and repelling directions with 
respect to a generic $1$-parameter subgroup in $\bA$. 

In these notes, we will see arrangements of both linear and 
affine rational hyperplanes and the components into which the
hyperplanes partition the real locus will play an important 
role. In both linear (or \emph{central}) and affine situations, 
these components are often 
called chambers or regions. To distinguish between 
the two, we will call the regions of central and affine 
arrangements \emph{cones} and \emph{alcoves}, respectively. 

We call a codimension $1$ stratum a \emph{wall}. A wall is 
a part of a hyperplane of the arrangement that separates two 
regions.

\subsubsection{}

The best definition of K\"ahler root of $X$ currently known to me is: 

\begin{Definition}
An effective curve class $\alpha$ is a positive K\"ahler root of $X$
if it appears in the sum \eqref{polesqmult}. 
\end{Definition}

It would be clearly desirable to have a more direct definition of a
root. One may try to define 
$$
\bigcup_{\alpha>0} \N \alpha\subset H_2(X,\Z)_\textup{eff}
$$
as those classes 
that remain effective in a codimension one deformation of $X$ 
as in Section \ref{s_deform}. There seems to be no simple way to 
pick out the roots among all of their multiples using this approach. 

Roots are also related to hyperplanes in $\Pic(X) \otimes_\Z \R$ along 
which stable envelope jump, see below and 
in particular \eqref{inclroot}. This approach requires a torus 
action which may not exist in general, and so cannot be 
used as a definition. 

\subsubsection{}

For $X=T^*G/B$, the components of the Steinberg 
variety 
$$
\St = \bigcup_{w\in W(G)}  L_w 
$$
are indexed by elements of the Weyl group of $G$. By the 
main result of \cite{BMO}, K\"ahler roots $\alpha$ of $X$ are the coroots
of $G_\textup{sc}$ with 
$$
L_\alpha = L_{s_\alpha} \,, 
$$
where $s_\alpha$ is the corresponding reflection. 

\subsubsection{}
Let $X=X_Q(\bv,\bw)$ be a Nakajima variety associated to a quiver $Q$ with 
vertex set $I$. 
The construction of \cite{MO}, the main points of which will be 
explained below, associates to $X$ a certain Lie algebra with 
a Cartan decomposition 
\begin{equation}
\fg_Q= \fh_Q \oplus \bigoplus_{\alpha} \left(\fg_Q\right)_\alpha\,, 
\label{gQcartan} 
\end{equation}
in which the root subspaces are finite-dimensional and are indexed
by roots $\alpha \in \Z^I$. Additionally
\begin{equation}
\left(\fg_Q\right)_{-\alpha} = \left(\fg_Q\right)_{\alpha}^\vee
\label{fgdual} 
\end{equation}
with respect to an invariant bilinear form on $\fg_Q$. 

For example, for the quiver $Q$ with one vertex and one loop 
one has $\fg_Q = \glh(1)$ and the roots of this algebra are 
all nonzero integers $\Z\setminus \{0\}$. 

As a corollary of the main result of \cite{MO}, we have 
\begin{equation}
 \textup{K\"ahler roots of $X_Q(\bv,\bw)$} = 
\textup{roots $\alpha$ of $\fg_Q$ such that $\pm\alpha\le \bv$} \,, 
\end{equation}
assuming $X_Q(\bv,\bw)\ne \varnothing$. Here $\alpha \le \bv$ means
that $\bv-\alpha \in \N^I$. 
In particular, 
$$
\textup{K\"ahler roots of $\Hilb(\C^2,n)$} = \{\pm1, 
\dots,\pm n\} \,.
$$

\subsubsection{}

In exceptional situations, it may turn out that the 
presentation of the form \eqref{polesqmult} is not unique because 
of special linear dependences that can exist between $L_\alpha$ and 
the identity operator. An example is 
\begin{equation}
\Hilb(\C^2,2) \cong T^*\bP^1 \times \C^2 \label{HilbTP} \,, 
\end{equation}
in which an attentive reader will notice that the two results just 
quoted give $\{\pm 1,\pm 2\}$ and $\{\pm 1\}$ as roots, respectively. 

Perhaps such nonuniqueness is always an indication 
of nonuniqueness of a lift to 
K-theory ? For example, the two points of view in \eqref{HilbTP} 
lead to \emph{different} lifts to K-theory, with different sets of 
singularities, matching the different sets of roots.

\subsubsection{}
There is a remarkable partial duality on equivariant 
symplectic resolutions which extends
$$
T^* G/B \xleftrightarrow{\quad \textup{Langlands} \quad 
} T^* {}^\textup{\tiny \sc L\,}\!G 
/ \,  {}^\textup{\tiny \sc L}\!B \,. 
$$
It has its origin in supersymmetric gauge theories, see 
e.g.\ \cite{BDGH,BLPW1} for a recent treatment, and is known in 
mathematics under the name of \emph{symplectic duality}. 
Among other things, it should exchange equivariant and 
K\"ahler roots. For example, the 
Hilbert schemes of points in $\C^2$ are 
self-dual.

\subsection{Braid groupoid}

\subsubsection{}

In classical Lie theory, a central role is played by the 
Weyl group of a root system. It is a natural group of 
symmetries which is large enough to act transitively on the 
set of cones. Further, the alcoves
of the \emph{affine arrangement}
\begin{equation}
\langle \alpha,x\rangle \in \Z \label{aff_arr}
\end{equation}
are permuted transitively by the 
affine Weyl group. 

There are natural finite groups of symmetries for both 
equivariant and K\"ahler roots. On the equivariant side, we 
have the Weyl group 
$$
\textup{Weyl}(\bA) = \textup{Normalizer}(\bA)/\bA 
$$
where the normalizer is taken inside $\Aut(X,\omega)$. 
On the K\"ahler side, there is an analog of the Weyl 
group constructed by Y.~Namikawa in \cite{NamWeyl}, see also e.g.\
the discussion of the topic in Section 2.2 of \cite{BPW}. 
In the affine case, one can take the semidirect product 
of this finite group with a suitable lattice. 

However, roots of either kind are intrinsically just not symmetric 
enough, e.g.\   for $X=\Hilb(\C^2,n)$, 
the alcoves are separated by the walls
$$
\left\{\frac ab\right\} \subset \R \,, 
\quad a\in \Z\,, b\in \{1,\dots,n\} 
$$ 
and they are very far from being transitive under $\{\pm \} \ltimes
\Z$.  

It is best to embrace the idea that 
different chambers and alcoves are different entities that 
only become 
related in some deeper way. 
 In this new paradigm, the place of 
a transitive symmetry group of a real hyperplane arrangement 
will be taken by the \emph{fundamental groupoid} of the 
same arrangement. I learned this point of view from 
R.~Bezrukavnikov.

\subsubsection{}

A groupoid is a category in which every morphism is invertible. 
The fundamental groupoid $\pi_1(U)$ of a topological space $U$ 
is formed by 
paths, taken up to homotopy that fixes endpoints.  Let 
$U=\C^n \setminus \bigcup H_i$ be the complement of 
a complexification of a real hyperplane arrangement. Here
$$
H_i = \{f_i(z)=0\}\,, 
$$
where $f_i$ is a real affine-linear function. The groupoid 
$\pi_1(U)$ was studied by Deligne \cite{Del1}, Salvetti \cite{Sal}, and many 
others. As in the survey \cite{Vas}, define 
$$
\Ct = \C_> \sqcup \C_< \sqcup \C_\uparrow
\sqcup \C_\downarrow 
$$
where 
$$
\C_{\gtrless} =\{\Re x  \gtrless 0 \} \,, \quad 
\C_{\uparrow\downarrow} =\{\Re x = 0\,, \Im x   \gtrless 0  \}  \,.
$$
For any sequence $*_1,*_2,\dots$ of symbols from $\{>,<,\uparrow,
\downarrow\}$, we set 
\begin{equation}
U_{*_1,*_2,\dots} = \bigcap f_i^{-1}(\C_{*_i}) \,.  \label{U**} 
\end{equation}
Nonempty intersections of the form \eqref{U**} decompose $U$. 
The Salvetti complex is, by definition, dual to this decomposition. 
It is a deformation retract of $U$.

We get 0-cells of the Salvetti complex when we choose 
$*_i\in \{\gtrless\}$ for all $i$; these corresponds to the alcoves
of $\R^n \setminus \bigcup H_i$ and give the objects of $\pi_1(U)$. 
The morphisms in $\pi_1(U)$ are generated by 1-cells of the 
Salvetti complex and those correspond to a wall between two 
alcoves plus a choice of over- or under-crossing, that is, a 
choice of $\uparrow$ versus $\downarrow$ for the corresponding 
equation $f_i$. 

\begin{figure}[!htbp]
  \centering
   \includegraphics[scale=0.4]{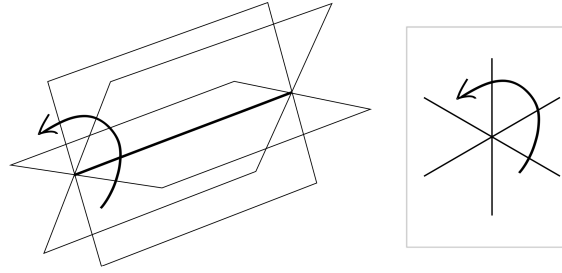}
 \caption{Product of morphisms around a stratum of codimension 2 
is the identity in the dynamical groupoid. In the fundamental 
groupoid, the imaginary parts of the  paths have to lie in one of 
the cones.}
  \label{f_walls_around}
\end{figure}

The relations correspond to 2-cells and 
those correspond an alcove $\nabla$ meeting several 
others along a stratum $S$ of codimension $2$, as in Figure 
\ref{f_walls_around}. 
We set 
$$
*_i = 
\begin{cases}
\gtrless\,, & f_i(S) \gtrless 0 \,, \\
\uparrow \downarrow\,, & f_i(S)=0, f_i(\nabla) \gtrless 0 \,. 
\end{cases}
$$
The corresponding $U_{**\dots}$ fibers over $S$ with fiber
$i \, \textup{Cone}_S\nabla$, where $\textup{Cone}_S\nabla$ 
denotes the tangent cone to $\nabla$ at $S$. 
We get a braid relations of the form shown in Figure \ref{f_walls_around}, 
where the path stays in $\R^n + i  \textup{Cone}_S\nabla$. 

\subsubsection{}
For future use, we note a generalization of $\pi_1(U)$ called 
the \emph{dynamical} groupoid of the arrangement. It has 
the same objects and for every wall $w$ between two alcoves
a single invertible matrix  
$B_w(f_i(z))$ which \emph{depends} on the variables $z\in \C^n$
through the equation $f_i$ of hyperplane $H_i$ containing $w$. 
The variables $z$ are called the dynamical variables. 

These matrices are required to satisfy one braid relation 
as in Figure \ref{f_walls_around} for any stratum of codimension 2. In other 
words, for the dynamical groupoid, there is a bijection 
\begin{equation*}
\begin{matrix}
\textup{objects} & \leftrightarrow & \textup{$\codim=0$ strata}\\
\textup{morphisms} & \leftrightarrow & \textup{$\codim=1$ strata} \\
\textup{relations} & \leftrightarrow & \textup{$\codim=2$ strata}
\end{matrix}
\end{equation*}
where strata refer to the stratification of $\R^n$. To go from 
the dynamical groupoid to the fundamental groupoid, we 
set
\begin{equation}
B_w^{\uparrow \downarrow} = \lim_{x\to \pm i\infty} B_w(x) \,,
\label{limB} 
\end{equation}
assuming this limit exists. 

We will meet the dynamical groupoids in exponential form, that is, 
as an arrangement of codimension $1$ subtori in an algebraic torus. 
The fundamental groupoid is then found from the values of 
$B_w(z)$ at the fixed points of a certain toric compactification. 

The Yang-Baxter equation \eqref{YB} is a very prominent example of 
a relation in a dynamical groupoid.

\subsubsection{}\label{s_K_arr} 
We now consider 
$$
U = H^2(X,\C) \setminus \,  \textup{K\"ahler arrangement} 
$$
where K\"ahler arrangement is the affine arrangement 
\eqref{aff_arr} 
corresponding to K\"ahler roots $\{\alpha\} \subset H_2(X,\Z)$. 
It is locally finite and periodic 
under the action of 
$$
\Pic(X) \cong H^2(X,\Z) \,.
$$
The fundamental groupoid $\pi_1(U)$ appears in two contexts, the 
relation between which is rather deep. These
are: 
\begin{itemize}
\item[---] monodromy of the quantum differential equation
  \eqref{qconn}, 
\item[---] autoequivalences of $\Db \Coh_{\bT} X$,
\end{itemize}
where $\bT \subset \Aut(X)$ is a maximal torus. 

{}In the early days of mirror symmetry, it was conjectured by 
M.~Kontsevich that any $X$ and $X'$ that have a common 
K\"ahler moduli space $\fK$ 
are derived equivalent and, moreover, there is 
a derived equivalence for any homotopy class of 
paths from $X$ to $X'$ in the regular locus of the quantum 
connection. This idea became more concrete with the 
advent of Brigeland's theory of stability conditions \cite{Br1,Br2}

More recently, it was conjectured by R.~Bezrukavnikov 
(with perhaps 
some infinitesimal input from the author of these notes, which is how 
is often attributed in the literature, see \cite{ABM,Et}) that, first, this works fully 
equivariantly\footnote{Equivariant mirror symmetry is a notoriously 
convoluted subject, see for example \cite{SeiSol}.} 
and, second, 
matches specific derived automorphisms 
of $X$ obtained via quantization in characteristic $p\gg 0$. 

\subsubsection{}\label{s_Iri} 

It is clear from definitions that the groupoid $\pi_1(U)$
acts by monodromy, that is, 
analytic continuation of solutions of \eqref{qconn}. 
Indeed, the natural isomorphism
\begin{equation}
H^2(X,\C)/H^2(X,\Z) \cong \Spec \Z H_2(X,\Z) \subset 
\fK \label{expK}
\end{equation}
takes $U$ to the regular locus of the quantum connection. 
Here $\Z H_2(X,\Z)$ is the group algebra of the lattice
$H_2(X,\Z)$. 

Let $z=0_X$ be the origin in the chart $\Spec \Z
H_2(X,\Z)_\textup{eff}$. Via the monodromy of \eqref{qconn}, 
the fundamental group of $U/\Pic(X)$ acts on the fiber at $0$. 
This fiber is $\Hd_\bT(X,\C)$, but we identify it with
$K_\bT(X)\otimes \C$ as in the work of Iritani \cite{Ir1}. 
Namely, one considers the map 
\begin{equation}
K_\bT(X) \owns \cF \to \boldsymbol{\Gamma} \cup 
(2\pi i)^{\frac{\deg(\,\cdot\,)-\dim_\C X}2} \chern \cF \in 
\Hd_\bT(X,\C)_\textup{mero} \label{Iri}
\end{equation}
where 
$$
\boldsymbol{\Gamma} = \prod_\textup{
$t_i$=Chern roots of 
$TX$}\Gamma(1+t_i)\,, 
$$
$\deg(\,\cdot\,)$ denotes the degree operator, and 
$\Hd_\bT(X,\C)_\textup{mero}$ denotes meromorphic functions on the 
spectrum of $\Hd_\bT(X,\C)$. Because of the identity 
$$
\Gamma(1+t) \, \Gamma(1-t) = \frac{2\pi i t} {e^{\pi i t} - e^{-\pi i
    t}}\,, 
$$
the map \eqref{Iri} is an analog of Mukai's map 
$\cF \to \chern(\cF) \sqrt{\textup{Td} X}$ in that it is, up to 
a scalar multiple, an isometry
for the natural \emph{sesquilinear} inner products in the source
and the target. It makes the monodromy act by 
unitary transformations of $K_\bT(X)\otimes \C$. 
Iritani's motivation was to find an \emph{integral} structure
in the quantum connection, and this will match nicely with what 
follows.

We note that the $q$-analog of the $\Gamma$-function 
is the function
\begin{align}
\Gamma_q(t) &= \prod_{n\ge 0} \frac1{1-q^n t}  
\notag \\
&=\tr \, \C[\textup{Maps}(\C(q)\to \C(t^{-1}))] 
\label{Gamma_q}
\end{align}
which is clearly relevant in enumerative K-theory of maps to $X$ 
because of the second line in 
\eqref{Gamma_q}. Here the equivariant 
weight of the source $\C(q)$ and of the target $\C(t^{-1})$ are
indicated in parentheses. See Section 8 of \cite{PCMI} for how 
these functions come up in the quantum difference equation 
for $X$.

\subsubsection{}
Quantization of equivariant symplectic resolutions is a very
fertile ground which is currently being explored by 
several teams of researchers, see for example 
\cite{ABM, BDMGN, BezF,BFG,BK1,BK2,BM,BezL,
BPW,BLPW2,Kaledin_q,Los1,Los2,Los3,MGN}.

A quantization of $X$ is, first, a sheaf 
$\cO_{\widehat{X}}$ of 
noncommutative algebras deforming the sheaf 
$(\cO_X, \{\,\cdot\,,\,\cdot\,\})$ of Poisson algebras
and, second, the algebra 
$$
\cX = \Gamma\left(\cO_{\widehat{X}}\right)
$$
of its global sections\footnote{The vanishing 
$H^i(\cO_X)=0$ for $i>0$ and any symplectic resolution 
implies the same for $\cO_{\widehat{X}}$, which is a very 
valuable property in the analysis of $\cX$.}. 
In this way, very nontrivial algebras can be constructed
as global section of sheaves of very standard algebras ---
a modern day replacement of the old dream of Gelfand 
and Kirillov \cite{GK}. 
Whenever $\cX$ is of finite homological 
dimension, one has a similar global description of  $\Db \cX
\textup{-mod}$. 

Quantizations $\cX_\lambda$ come in families parametrized by 
the same data $$
\lambda \in \Pic(X) \otimes \textup{ground field}
$$
as the commutative deformations. 
For Nakajima varieties, these can be described explicitly as 
so-called quantum Hamiltonian reductions \cite{EtCal}. 

\subsubsection{}

While these notes are certainly not the place to survey the 
many successes of quantization of symplectic resolutions, 
certain important features of quantizations in characteristic 
$p\gg 0$, developed by Bezrukavnikov and Kaledin, 
 are directly related to our narrative. 

For an integral quantization parameter $\lambda$ away from certain
walls, that is for  
$$
\lambda \in \Pic(X) \cap p U'
$$
where 
$$
U' = \Pic(X) \otimes \R \, \setminus \,\, \textup{neighborhood
 of a periodic 
hyperplane arrangement} \,,
$$
the theory produces an equivalence 
\begin{equation}
\Db \Coh_{\bA^{(1)},\textup{supp}} X^{(1)}  
\xleftrightarrow{\quad\sim\quad} \Db
\cX_\lambda\textup{-mod}_{\bA,\textup{supp}}
\label{Dequiv} 
\end{equation}
where 
$$
\bA\subset\Aut(X,\omega)=\Aut(\cO_{\widehat{X}_\lambda})
$$
is a maximal torus, $X^{(1)}$ denotes Frobenius twist, and 
one needs to impose a certain condition on the supports of 
sheaves and modules in \eqref{Dequiv}. We will assume that
the projection of the set-theoretic support to the affinization $X_0$ 
is contracted to $0\in X_0$ by a certain $1$-parameter
subgroup in $\bA$; this means a version of the 
category $\cO$ on the quantization side. 

On the commutative side, we have a larger torus $\bT \supset \bA$ 
of automorphisms, thus the quantization side gets an extra grading
so that \eqref{Dequiv} is promoted to a $\bT$-equivariant statement. 
The extra grading  on $\cX_\lambda\textup{-mod}$
captures quite subtle 
representation-theoretic information. In fact, such graded 
lifts of representation categories are absolutely central to a lot of 
progress in modern 
representation theory, see e.g.\ \cite{BGS,Soe}. 

\subsubsection{}\label{s_Bez} 

The map \eqref{expK} send $U'\subset \Pic(X)\otimes \R$ to the 
central torus $|z|=1$ in $\fK$ and one can associate the 
equivalence \eqref{Dequiv} to a straight 
path from the point $z=0_X\in \fK$ 
to a component of $U'$. 

One, meaning Bezrukavnikov, conjectures that:
\begin{enumerate}
\item[---] the hyperplane arrangement in the 
definition of $U'$ is the affine K\"ahler  
arrangement\footnote{For uniformity, it is convenient to 
allow walls with trivial monodromy and derived equivalences, 
respectively. For example, the point $z=1$ is not a singularity 
for the Hilbert scheme of points, as the pole there is cancelled 
by the scalar operator in \eqref{polesqmult}.}
 as in 
Section \ref{s_K_arr}, 
\item[---] the equivalences \eqref{Dequiv} define a 
representation of $\pi_1(U)$, and
\item[---] their action on $K_\bT(X)$ equals the monodromy of 
the quantum differential equation, 
\end{enumerate}
where we use Iritani's map as in Section \ref{s_Iri} to 
lift monodromy to K-theory. 

The second part of this conjecture
means the following. For $\lambda \in p U'$, the algebras
$\cX_\lambda$ are canonically identified for all flops of $X$ and, 
for two quantization parameters $\lambda, \lambda'$ separated by a wall 
$w$, we can consider the induced equivalence, that is, the 
horizontal row in the following diagram: 
\begin{equation}
\label{triangDb} 
\xymatrix{
& \Db \Coh  X^{(1)}  \ar[dr] \ar[dl] \\
\Db \cX_\lambda\textup{-mod} \ar[rr]^{B_w^{\uparrow\downarrow}}
&& 
\Db \cX_{\lambda'}\textup{-mod}  \,. 
}
\end{equation}
The corresponding paths fall into two homotopy classes, 
according to the equation of the wall being effective or 
minus effective in $X$. Therefore, one expects all $X$ in 
the same class to induces the same equivalence
$B_w^{\uparrow}$, 
resp.\ $B_w^{\downarrow}$, in the diagram \eqref{triangDb}. 
For recent results in this direction, see \cite{Boger}. 

Note that shifts $\lambda$ by the lattice $p\Pic(X)$ amount to 
twists of $\Db \Coh  X^{(1)}$ by a line bundle.  
This gives very interesting factorizations of twists by 
line bundles in the group $\Aut \Db \Coh X$. 

\subsubsection{}

For representation theory of $\cX$, the monodromy group 
plays the role of Hecke algebra in the classical Kazhdan-Lusztig
theory. It packages very valuable 
representation-theoretic information, which e.g.\ includes the 
classification of irreducible $\cX$-modules according 
to their size. The latter problem inspired further conjectures by 
Etingof \cite{Et}, which were proven in many important 
cases \cite{BezL,Los1,ShanV}. The general statement about monodromy 
is currently known \cite{BO} for Nakajima varieties $X$ such that 
$\dim X^\bA=0$. 

As we will see, the monodromy of the quantum differential 
equation \eqref{qconn} 
lies somewhere between the differential equation 
itself and its $q$-difference analog. Possible categorical 
interpretations of the quantum $q$-difference equations will be 
discussed below. 

\section{Stable envelopes and quantum groups}
\label{s_st} 

\subsection{Stable envelopes} 

\subsubsection{}

Let $X(\bw)$ be a Nakajima variety as in \eqref{Xw}. 
Our goal in this section is to produce interesting correspondences
in $X(\bw) \times X(\bw)$ which will result in actions of quantum 
groups and will be used to describe solutions to enumerative
problems. 

Note that $X(\bw)$ is disconnected, and so a correspondence 
in it is really a collection of correspondences in 
$X(\bw,\bv) \times X(\bw,\bv')$. Even if one is interested
in correspondences in $X \times X$ for a 
given symplectic resolution $X$, it 
proves beneficial to use correspondences $X \times X'$ as 
ingredients in the construction. In fact, perhaps the main 
obstacle in generalizing what is known for Nakajima varieties
to general symplectic resolutions $X$ is the shortage of 
natural relatives $\{X',X'',\dots\}$ with which $X$ could meaningfully 
interact. 

A rare general construction that one can use is the following. 
Let 
$$
\bA \subset \Aut(X,\omega)
$$
be a torus, not necessarily maximal. Then a choice of a $1$-parameter
subgroup $\sigma: \Ct \to \bA$ determines a 
Lagrangian submanifold 
\begin{equation}
\Attr =\left\{(x,y),\,  \lim_{t\to 0} \sigma(t) \cdot x = y 
\right\} \subset X \times
X^\bA \,. \label{Attr} 
\end{equation}
This is a very familiar concept, examples of which include e.g.\ 
conormals to Schubert cells in $T^*G/P$ or the loci of extensions in 
the moduli of framed sheaves. 

The choice of $\sigma$ matters only 
up to the cone $\fC$ 
$$
d\sigma \in \fC  \subset \Lie \bA 
$$
cut out by the equivariant roots of $X$ and containing $d\sigma$. 
We already saw these cones in Section \ref{scones}. 
We will write $\Attr_\fC$ when we need to stress this dependence. 

\subsubsection{}\label{sStab_spec} 

The submanifold $\Attr$ cannot be used for our purposes for the simple
reason of not being closed. Using it closure, and especially the 
structure sheaf of the closure, runs into all the usual problems with 
closure in algebraic geometry. These include $\overline{\Attr}$ 
not being stable against perturbations, that is, not fitting into a 
family for $\bA$-equivariant deformations $X_\lambda$ of $X$. 

For generic $\lambda\in H^2(X,\C)$, the deformation $X_\lambda$ 
is affine and \eqref{Attr} is closed for it. In cohomology, one 
gets good results by closing this cycle in the whole family, see
Section 3.7 in \cite{MO}. In the central fiber $X$, this 
gives a Lagrangian cycle 
$$
\Stab_\fC \subset X \times X^\bA \,,
$$
called \emph{stable envelope}, 
supported on the full attracting set $\Attr^f$. By definition 
$(x,y)\in \Attr^f$ if $x$ can be joined to $y$ by a chain of 
closures of attracting $\bA$-orbits. With hindsight, one 
can recognize instances of this construction in such 
classical works as \cite{Ginz}. 

In practice, it is much more useful to have a characterization of 
$\Stab_\fC$ in terms that refer to $X$ alone instead of 
perturbations and closures. This becomes crucial in 
$\bT$-equivariant K-theory, where $\bT\supset \bA$ is a 
torus which scales $\omega$ nontrivially. 
Since there are no $\bT$-equivariant deformations of $X$, 
a perturbation argument does not yield a well-defined 
$\bT$-equivariant K-class. In fact, as we will see, 
stable envelopes in K-theory crucially depend on certain 
data other than just a choice of a cone $\fC$. 

\subsubsection{}
The fixed locus has many connected components $X^\bA = 
\bigsqcup F_i$ and, by definition, 
$$
\left(\supp \Attr^f_\fC\right)^\bA = 
\bigsqcup_{i\ge j} F_i \times F_j \,,
$$
where $i\ge j$ refers to the partial order on components of 
$X^\bA$ determined by the relation of being attracted by $\sigma$. 
The main idea in definition of stable envelopes both in 
cohomology and $K$-theory is to require
that 
\begin{equation}
\Stab\Big|_{F_i \times F_j} \quad \textup{is smaller than} 
\quad \Stab\Big|_{F_i \times F_i}\,, \label{StabStab} 
\end{equation}
for $i>j$ in a sense that will be made precise presently. 

The torus $\bA$ doesn't act on the fixed locus, so the 
restriction in \eqref{StabStab} in a polynomial in equivariant 
variables with values in $\Hd_{\bT/\bA}(F_i \times F_j)$ and 
$K_{\bT/\bA}(F_i \times F_j)$, respectively. We denote 
$\deg_\bA$ the degree of this polynomial, defined as follows. 
In cohomology, we have a polynomial on $\Lie \bA$, with the 
usual notion of degree. We require
\begin{equation}
\deg_\bA \Stab\Big|_{F_i \times F_j} \quad < 
\quad \deg_\bA \Stab\Big|_{F_i \times F_i}
= \frac12 \codim F_i \,, \label{StabStab2} 
\end{equation}
and a simple argument, see Section 3 in \cite{MO},
 shows that an $\bA$-invariant 
cycle $\Stab_\fC$ supported on 
$\Attr^f$, equal to $\Attr$ near diagonal, 
and satisfying the degree bound \eqref{StabStab2} 
\begin{enumerate}
\item[---] is always unique; 
\item[---] exists for rather general algebraic 
symplectic varieties $X$, in particular, for 
all symplectic resolutions; 
\item[---] for symplectic resolutions, may be obtained by a 
specialization argument from Section \ref{sStab_spec}. 
\end{enumerate}
It can be constructed inductively by a version of 
Gram-Schmidt process, the left-hand side in \eqref{StabStab} is then 
interpreted as the remainder of division by the right-hand side. 
Since we are dealing with multivariate polynomials, this is 
somewhat nontrivial. The key observation here is that 
for $\bA \subset Sp(2n)$, the class of any $\bA$-invariant 
Lagrangian $L\subset
\C^{2n}$ in $H^{2n}_\bA(\C^{2n})$ is a $\Z$-multiple of the class of any 
fixed linear Lagrangian.  A naive argument like this doesn't work 
in $K_\bA(\C^{2n})$ and things become more constrained there.

\subsubsection{}
In K-theory, we deal with polynomials on $\bA$ itself, for which 
the right notion of degree is given by the Newton polygon 
$$
\deg_\bA \sum c_\mu a^\mu = 
\textup{Convex hull} \left( \{\mu, \, c_\mu \ne 0\} \right) 
\subset \bA^\wedge \otimes_\Z \R 
$$
considered \emph{up to translation}. The natural ordering 
on Newton polygons is that by inclusion and to allow for the up-to-translation ambiguity we require stable envelopes in 
K-theory to satisfy the condition 
\begin{equation}
\deg_\bA \Stab\Big|_{F_i \times F_j} \quad \subset 
\quad \deg_\bA \Stab\Big|_{F_i \times F_i}
+ \textup{shift}_{ij} \label{StabStab3} 
\end{equation}
for a certain collection of shifts 
$$
\textup{shift}_{ij} \in \bA^\wedge \otimes_\Z \R \,.
$$
Condition \eqref{StabStab3} is known as the \emph{window}
condition. 
The uniqueness of K-classes satisfying it is, again, immediate. However,
existence, in general, is not guaranteed. If the rank of $\bA$ is one, 
one can argue inductively using the usual division with remainder. 
However, for application we have in mind it is crucial to allow for 
tori of rank $>1$, and for them existence can only be shown for 
special shifts associated to fractional line bundles on $X$. 

Let $\cL \in \Pic_\bA(X) \otimes_\Z \R$ be a fractional
$\bA$-linearized
line bundle. Its restriction $\cL\Big|_{F_i}$ to a component of the fixed 
locus has a well-defined weight in $\bA^\wedge \otimes_\Z \R$ 
and we define 
$$
\textup{shift}_{ij} = \textup{weight} \, \cL\Big|_{F_i} - 
\textup{weight} \, \cL\Big|_{F_j} \,.
$$
Note that the choice of linearization 
cancels out and the shift depends on just a 
fractional line bundle $\cL \in \Pic(X) \otimes_\Z
\R$, called the \emph{slope} of the stable envelope. 

There is some fine print in the correct normalization of 
K-theoretic stable envelopes at the diagonal. It involves the notion 
of a \emph{polarization}, see \cite{PCMI} for details. In fact, even in 
cohomology, it is best to require $\Stab = \pm \Attr$ near 
the diagonal, depending on a polarization, see \cite{MO}. 
With such normalization and for a certain class of $X$ that 
includes all symplectic resolutions,  we have 

\begin{Theorem}[\cite{MO2,AO,HLMO}] 
There exists unique stable envelope 
\begin{equation}
\Stab_{\fC,\cL} \in K_\bT(X \times X^\bA) \label{StabCL} 
\end{equation}
for any cone $\fC$ and any slope $\cL$
 away from a certain locally finite 
$\Pic(X)$-periodic rational 
hyperplane arrangement in $\Pic(X)\otimes_\Z\R$. 
\end{Theorem}

Proofs in the literature deduce this from 
more general statements. For Nakajima varieties, 
the existence of more general \emph{elliptic} 
stable envelopes is shown in \cite{AO} , see 
Section 3.8 there for specialization to K-theory. 
The paper 
\cite{HLMO} gives a 
Gram-Schmidt-style existence proof of \emph{categorical} 
stable envelopes.

Note that for generic shifts, the inclusion in \eqref{StabStab3} 
is necessarily strict as the inclusion of an integral polytope into 
a very nonintegral one. Therefore, the $\cL$-dependence of 
$\Stab_{\fC,\cL}$ is locally constant and can only jump if 
a lattice point gets on the boundary of the polytope on the right 
in \eqref{StabStab3}. 

There is a direct relation between stable envelopes and the 
monodromy of \eqref{qconn} proven in \cite{AO} in the 
more general $q$-difference case. It implies, in particular 
that 
\begin{equation}
\left\{
 \begin{matrix}
\textup{hyperplanes in $\Pic(X)\otimes_\Z\R$ } \\
\textup{where $\Stab_{\fC,\cL}$ jump}
\end{matrix}
\right\} 
\, {\Huge \subset}  \, 
\left\{
\begin{matrix}
 \textup{K\"ahler arrangement} \\
\textup{for $X$} 
\end{matrix} 
\right\}\,, \label{inclroot} 
\end{equation}
and I think that we have an equality in \eqref{inclroot} aside from 
trivialities like $X=X_1 \times X_2$ with $\bA$ acting on $X_2$ only. 
Also note that the left-hand side in \eqref{inclroot} does not 
depend on $\fC$. 

\subsubsection{}
Piecewise constant dependence on a fractional line bundle should 
certainly activate neurons in areas of the cerebral cortex responsible 
for multiplier ideal sheaves and related objects. Indeed, the
similarity\footnote{It is easy to suspect an actual link here, and 
I would be grateful to experts for explaining what it is.} 
 goes further in that 
\begin{enumerate}
\item[---] window conditions are 
preserved by proper push-forwards; 
\item[---]  in a normal-crossing 
situation (e.g.\ for hypertoric varieties), stable envelope is a 
 twist of a structure sheaf 
by a rounding of $\cL$. 
\end{enumerate}
In particular, Nakajima varieties are related by certain 
Lagrangian \emph{abelianization} correspondences to 
smooth hypertoric varieties, see \cite{HP}, and this has been used
to compute stable envelopes starting in \cite{Shen,Sm2} and 
all the way to the elliptic level of generality in \cite{AO}. 

\subsubsection{}

Favorable properties of stable envelopes are shared by 
their categorical lifts, which are defined as functors
\begin{equation}
\mathbf{Stab}_{\fC,\cL} : \Db \Coh_\bT X^\bA \to 
 \Db \Coh_\bT X \label{stab_cat}
\end{equation}
satisfying the same window conditions for the derived 
restriction to $X^\bA$. For $\rk \bA =1$, a detailed 
study of these, with many applications, can be found the 
in the work of Halpern-Leistner \cite{HL}, see also \cite{BFK}. 
In this case, one doesn't need the shift to come from 
a line bundle. The case of tori of higher rank is considered
in \cite{HLMO}. 

\subsection{R-matrices}

\subsubsection{}
If a construction in algebraic geometry requires additional 
choices, such as the choice of a cone $\fC$ and a slope $\cL$ 
in the definition of stable envelopes, then this can be 
viewed as a liability or as an asset. The latter point of view
is gaining popularity as people find more and more applications 
of wall-crossing phenomena of various kind. 
In the case at hand, it will be beneficial to study how stable 
envelopes change as the we cross the wall from 
$\fC$ to another chamber $\fC'\subset\Lie \bA$, or as we 
cross a wall of the K\"ahler arrangement from $\cL$ to 
$\cL'$. 

To formalize this, we introduce
\begin{equation}
R^{\cL}_{\fC'\leftarrow \fC} = 
\Stab_{\fC',\cL}^{-1} \, \circ \, \Stab_{\fC,\cL} \label{RC} 
\end{equation}
and similarly 
\begin{equation}
R^{\fC}_{\cL'\leftarrow \cL} = 
\Stab_{\fC,\cL'}^{-1} \, \circ \, \Stab_{\fC,\cL} \,. \label{RL} 
\end{equation}
Stable envelopes are 
isomorphisms\footnote{The inverse of a stable envelope is a 
transpose of another stable envelope, see Section 9 in \cite{PCMI}.} 
after localization, so 
$$
\textup{all $R$-matrices} \subset \End(K_\bT(X^\bA)) \otimes 
\Q(\bT) \,.
$$
They \emph{depend} on equivariant variables $a\in \bA$ for the 
torus $\bA$ which does not act on $X^\bA$, and this is, again, a
bonus.  We can expand R-matrices as $a\to 0$, where $0$ is a
fixed points of a toric compactification $\overline{\bA} \supset \bA$ 
that correspond to the fan $\{\fC_i\}$ of chambers. The 
coefficients of this expansion give countably many operators 
acting on \emph{integral} $\bT/\bA$-equivariant K-theory of 
$X^\bA$, that is, we can 
 view $R$-matrices as interesting generating functions 
for correspondences in $X^\bA$. 

\subsubsection{}

Let $\fC$ and $\fC'$ be separated by a wall in $\Lie \bA$. This 
wall is of the form $\Lie \bA'$ for some codimension 1 subtorus
$\bA' \subset \bA$. One shows, see Section 9.2 in \cite{PCMI} that 
\begin{equation}
R_{\fC' \leftarrow \fC} \textup{  for $\bA$-action on $X$  } 
= \, 
R_{\fC' \leftarrow \fC} \textup{  for $\bA/\bA'$-action on $X^{\bA'}$
} 
\label{RRAA} 
\end{equation}
and, in particular, $R_{\fC' \leftarrow \fC}$ factors through the 
map $\bA \to \bA/\bA'$, that is, through the equation of the wall. 
This gives

\begin{Proposition}
Matrices $R_{\fC' \leftarrow \fC}$ form a representation of the
dynamical groupoid for the central arrangement defined by 
the equivariant roots of $X$. 
\end{Proposition}

Here central is used as opposite of affine: we equate equivariant 
roots to $0$ and not to an arbitrary integer. 

\subsubsection{}

In the situation of Section \ref{stensNak}, the root hyperplane
are $a_i/a_j=1$ and the groupoid relation becomes
the Yang-Baxter equation. 

This geometric solution of the YB 
equation produces a geometric action of a quantum group
 $\cA = \cU_\hbar(\gh)$ for 
a certain Lie algebra $\fg$ as indicated in Section 
\ref{s_introR} . See \cite{MO,OS} for 
the details of the reconstruction.

The parallel procedure in cohomology constructs a certain 
degeneration $\bY(\fg)$ of $\cU_\hbar(\gh)$ known as the Yangian 
of $\fg$. This is a Hopf algebra deformation of $\cU(g[t])$ and 
the universal enveloping of $\fg \subset \fg[t]$ 
\begin{equation}
\cU(\fg) \subset \bY(\fg) \label{UgY} 
\end{equation}
remains undeformed as a Hopf algebra. 

\subsubsection{}
The structure of $\cA$ may be analyzed using the following 
useful connection between R-matrices of the two kinds 
\eqref{RC} and \eqref{RL}.  

We can include $\cL$ into an 
infinite sequence of slopes 
\begin{equation}
\dots, \cL_{-2}, \cL_{-1}, \cL_0 = \cL, \cL_1, \cL_2, \dots 
\label{path} 
\end{equation}
that go from a point $\cL_{-\infty}$ on the minus ample 
infinity of $\Pic(X)$ to a point $\cL_{\infty}$ on the ample infinity
of $\Pic(X)$. A choice of such path is analogous to a choice 
of the factorization of the longest element in Weyl group in the 
usual theory. 

Tautologically, 
\begin{equation}
R^{\cL}_{(-\fC)\leftarrow\fC} = 
\overrightarrow{\prod_{i < 0 }} R^{-\fC}_{\cL_{i+1} \leftarrow \cL_i}
\,\,   R_{(-\fC,\cL_{-\infty})\leftarrow (\fC,\cL_\infty)} 
\,\, \overleftarrow{\prod_{i \ge 0 }}
R^{\fC}_{\cL_{i+1} \leftarrow \cL_i} \,,
\label{pathfact} 
\end{equation}
assuming the term in the middle makes sense, that is, infinite
product converge in a suitable topology. This is indeed the case: 

\begin{Proposition}[\cite{OS}] 
 The operator $R_{(-\fC,\cL_{-\infty})\leftarrow (\fC,\cL_\infty)}$ is
 well-defined and acts by an operator of multiplication 
in $K_\bT(X^\bA)$. 
\end{Proposition}

In fact, this operator acts by 
a very specific Schur functor of the normal bundle $N_{X/X^\bA}$ 
to the fixed locus. 

Formula \eqref{pathfact} 
 generalizes factorization of R-matrices constructed 
in the  conventional 
theory of quantum groups, see for example
\cite{KirResh,KhorTol, LevSoi}. 

Recall that $\cU_\hbar(\gh)$ is spanned by coefficients in
$u\in\bA\cong \Ct$ 
of the matrix coefficients of R-matrices
\eqref{RC}. In this procedure, 
one can, in fact, pick out the matrix coefficients of each 
individual term in the factorization \eqref{pathfact} and 
this gives 
\begin{equation}
\cU_\hbar(\fg_w) \hookrightarrow \cU_\hbar(\fgh)
\label{Ugw} 
\end{equation}
for each wall $w$ crossed by the path \eqref{path}. These deform 
the decomposition 
\begin{equation}
\fgh = \fg[t^{\pm 1}] = \fhh \oplus \bigoplus_{w} \fg_w \,, 
\quad 
\fg_w = \bigoplus_\textup{$\alpha=-n$ on $w$}  \fg_\alpha t^n \,,
\label{fgdecomp} 
\end{equation}
in which the $\fhh$ part corresponds to the middle factor in 
\eqref{pathfact}. Note that the subalgebra \eqref{Ugw} 
depends on 
\begin{enumerate}
\item[---] the wall $w$ between two alcoves and not only on the 
hyperplane containing it and, moreover, on 
\item[---] the choice of the path \eqref{path}. 
\end{enumerate}
These features are already familiar from the classical theory of 
quantum groups. Concretely, in textbooks, quantum groups are 
usually presented by generators and relations and, consequently,
come equipped with root subalgebras corresponding to simple roots. 
Constructing other root subalgebras requires choices just like 
we had to make in \eqref{path}. 

This argument also shows that the quantum group constructed
using the R-matrix \eqref{pathfact} is independent of the 
choice of $\cL$. 

\subsection{Enumerative operators}

\subsubsection{}
We go back to curve-counting in cohomology and reexamine 
the operator \eqref{qmult} for Nakajima varieties. 
Recall the inclusion \eqref{UgY}. The root decomposition 
in \eqref{gQcartan} is according to the components of 
$X(\bw)$, that is, 
\begin{equation}
\fg_\alpha: \Hd(X(\bv,\bw)) \to \Hd(X(\bv-\alpha,\bw))\,,
\label{actga} 
\end{equation}
and all these root spaces are finite-dimensional. In parallel 
to Kac-Moody theory, one shows that all roots are either 
positive or negative,  that is, 
$$
\textup{roots of $\fg$} \subset  \N^I \cup (-\N)^I \,.
$$
Duality \eqref{fgdual} gives a canonical element 
$$
\Casimir_\alpha \in \,:\!\fg_{-\alpha} \fg_\alpha\!:\, \subset \cU(\fg) 
$$
which preserves $\Hd X(\bv,\bw)$ and annihilates it unless 
$|\alpha| \le \bv$ because of \eqref{actga}. Normal ordering 
here means that we act by lowering operators first. 

The main result of \cite{MO} is the following 

\begin{Theorem}
\label{tH}
For Nakajima quiver varieties, formula \eqref{polesqmult} holds with
$$
L_\alpha = - \Casimir_\alpha  \,.
$$
\end{Theorem}

In particular, this gives a new computations of the quantum 
cohomology of Hilbert schemes of ADE surfaces, first obtained in 
\cite{OPhilb,MOb} by a more direct analysis. 

\subsubsection{}

The proof of Theorem \ref{tH} given in \cite{MO} is rather indirect. 
Recall that the group $\prod GL(W_i)$ in \eqref{GLGL} 
acts on $X(\bv,\bw)$,
and this action is nontrivial if $\sum \bw_i > 1$. There
is a flat difference connection in the corresponding 
equivariant variables which commutes with the quantum 
connection \eqref{qconn}. It is constructed geometrically 
counting sections in certain twisted $X$-bundles over $\bP^1$. 
Operators of this kind, often called \emph{shift operators}, find 
many other applications in enumerative geometry, see for 
example \cite{Ir2}. 

The main step of the proof identifies this commuting connection 
with what is known as the \emph{quantum Knizhnik-Zamolodchikov} 
equations. These were introduced by Frenkel and Reshetikhin in 
\cite{FrenkelReshetikhin} in the context of 
 quantum loop algebras of finite-dimensional Lie algebras. 
They make sense for any solution of the Yang-Baxter equation with a
spectral parameter and, in particular, for our geometrically 
constructed R-matrices.
The support and the degree conditions 
satisfied by stable envelopes enter the argument at this step. 

The K-theoretic 
computation in \cite{PCMI,OS}, while much more involved on 
the technical level, follows the same general strategy: one 
identifies a commuting difference connection first, and uses 
that to constrain the connection of principal geometric 
interest.

\subsubsection{}

One may phrase the answer in Theorem \ref{tH} as 
the identification of the quantum connection \eqref{qconn} with 
the trigonometric \emph{Casimir connection} for 
the Yangian $\bY(\fg)$. For finite-dimensional Lie algebras, 
such connections were studied by Toledano Laredo in \cite{TL}, see also
the work \cite{TV1,TV2} by Tarasov and Varchenko.

A connection with Yangians was expected from the 
ideas of Nekrasov and Shatashvili. The precise identification 
with the Casimir connection was suggested by Bezrukavnikov, 
Etingof, and their collaborators. Back then, it was 
further predicted by Etingof that the
correct K-theoretic version of the quantum connection 
should be a generalization of the dynamical Weyl group 
studied by Tarasov, Varchenko,
Etingof, and others for finite-dimensional
Lie algebras $\fg$, see e.g.\ \cite{TV3,EV}. 
 Indeed, for such $\fg$, 
Balagovic showed \cite{Bal}
 that the dynamical $q$-difference equations
degenerate to the Casimir connection in an appropriate
limit. 

At that moment in time, there was neither 
geometric, nor a representation-theoretic construction 
of the required $q$-difference connection. 

\subsubsection{}

Referring the reader to \cite{CKM} for all details, we will define 
a quasimap 
$$
f: C \dashrightarrow  X = \widetilde{X} \rdd G
$$
to a GIT quotient as a section of a bundle of prequotients 
$\cP\times_G X$, up to isomorphism. Here $\cP$ is a 
principal $G$-bundle over $C$ which is a part of the data and 
is allowed to vary. Stability conditions for quasimaps are 
derived from stability on $ \widetilde{X}$; we will use the
simplest one that requires the value of the quasimap be 
stable at the generic point of $C$. 

To have an evaluation map at a specific point $p\in C$ one introduces 
quasimaps relative $p$, formed by diagrams of the form 
\begin{equation}
\xymatrix{
p'\ar@{|->}[r]\ar@{|->}[d]
&C' \ar@{-->}[rr]^f \ar[d]_\pi  && \widetilde{X} \rdd G \,, \\
p\ar@{|->}[r] &C
}\label{relquasi} 
\end{equation}
in which $\pi$ is an isomorphism over $C\setminus \, p$ and 
contracts a chain of rational curves joining 
$p'$ to $\pi^{-1}(C\setminus \, p)$. The evaluation map 
records the value at $p'$. 

For quasimaps to $X(\bv,\bw)$,  we have  $G=\prod GL(V_i)$ and 
so a principal $G$-bundle is the same as a collection of 
vector bundles of ranks $\bv$.  For a regular map $f$, these 
are the pullbacks $f^* \cV_i$ of the tautological vector bundles
$\cV_i$ 
on $X$. We denote them by $f^* \cV_i$ for all quasimaps; they 
are a part of the data of a quasimap. 

Consider the line bundle 
$$
\cL_i = \det \cV_i \in \Pic(X) \,. 
$$
For our difference connection, we need a correspondence in $X$ 
that generalizes the operator $\otimes \cL_i$ in $K_\bT(X)$. 
We take
\begin{equation}
\bfM_{\cL_i} = \sum_d z^d \, 
\ev_* \left( \tO_{\cM_d} \otimes \det \Hd( f^*\cV_i \otimes \pi^*\cO_{p_1})
\right) \, \bfG^{-1} \,,  \label{defbfM} 
\end{equation}
where $\cM_d$ is the moduli space of quasimaps 
$$
(\bP^1,0,\infty) \cong (C,p_1,p_2) \dasharrow \widetilde{X} \rdd G
$$
of degree $d$ and $\bfG$ is the glue matrix from \eqref{defglue}. 
See \cite{PCMI} for detailed discussion of how one 
arrives at the definition \eqref{defbfM} and what is the 
geometric significance of the solutions of the corresponding 
$q$-difference equations. 

\subsubsection{}

On the algebraic side, recall the subalgebras \eqref{Ugw} 
associated to a wall of the K\"ahler arrangement. 
It is by itself a quantum group, in particular, it has 
its own R-matrix $R_w$, which is the corresponding 
factor in \eqref{pathfact} with a suitable normalization. 
This is not a quantum loop algebra, in fact it is a quantum 
group of rank 1. Correspondingly, $R_w$ does not have a 
spectral parameter and satisfies the YB equation \eqref{YB} 
for constant (in $u$) operators. Examples of these 
wall subalgebras are $\cU_\hbar(\mathfrak{sl}(2))$ and
the quantum Heisenberg algebra. 

By definition, quantum Heisenberg algebra has two 
group-like invertible elements 
$$
\Delta H = H \otimes H \,, \quad \Delta K = K \otimes K \,,
$$
of which $K$ is central, while 
$$
H E = \hbar \, E H \,, \quad H F = \hbar^{-1} F H\,. 
$$
Here the generators $E$ and $F$ satisfy
$$
\left[ E, F \right] = \frac{K-K^{-1}}{\hbar-\hbar^{-1}} 
$$
and comultiply as follows
\begin{align*}
\Delta E & = E \otimes 1 + K^{-1} \otimes E  \,. \\
\Delta F & = F \otimes K + 1 \otimes F \,.
\end{align*}
This algebra acts on $\Z[x,\hbar^{\pm 1}]$ by 
$$
E = x \,, \quad F = - \frac{d}{dx}\,, \quad H = 
\hbar^{x \tfrac{d}{dx}}\,, \quad 
K = \hbar \,. 
$$
Its R-matrix is 
$$
R = \hbar^{\Omega} \, \exp\left( - (\hbar - \hbar^{-1}) F \otimes E
\right) \,,
$$
where 
$$
\Omega = - \log_\hbar H \otimes \log_\hbar H  = - x \tfrac{d}{dx}
\otimes x \tfrac{d}{dx} \,. 
$$
In its geometric origin, the $\Z$-valued operator $\Omega$ records
the codimension of the fixed locus. 

For quivers of affine type, every root is either real, in 
which case the wall subalgebra is 
isomorphic to  $\cU_\hbar(\mathfrak{sl}(2))$, or 
has zero norm. In the latter case, all multiples of the 
root enter $\fg_w$ in \eqref{fgdecomp}
 and $\cU_\hbar(\fg_w)$ is 
a countable product of Heisenberg subalgebras, with their tori 
identified.  For example, all root subalgebra of the algebra 
$\cU_\hbar(\widehat{\widehat{\mathfrak{gl}}}(1))$ associated
with the quiver with one vertex and one loop are of this type. 
The structure of this algebra has been studied by many authors, 
see for example \cite{BurSch,FeiJ1,FeiJ2,Neg}. 

\subsubsection{}
With any R-matrix, we can associate 
\begin{equation}
\textup{qKZ operator} = ( z^{\log_\hbar H} \otimes 1 ) \, R = 
( z^{\log_\hbar H} \otimes 1 ) \, \hbar^{\Omega} + \dots 
\label{qKZoper} 
\end{equation}
acting in the tensor product of two $\cU_\hbar(\fg_w)$-modules. 
Here dots stand for an upper-triangular operator. 

Since $\cU_\hbar(\fg_w)$ is not a loop algebra, this is not a
difference operator. Instead of a fundamental solution of 
a difference equation, we can ask for a unipotent matrix $\bJ$ 
that conjugates \eqref{qKZoper} to its eigenvalues 
\begin{equation}
( z^{\log_\hbar H} \otimes 1 ) \, R \, \bJ = 
\bJ \, ( z^{\log_\hbar H} \otimes 1 ) \, \hbar^{\Omega} \,.
\label{Jconj} 
\end{equation}
Such operator exist universally and for the quantum Heisenberg 
algebra it equals 
$$
\bJ = \exp\left( - (\hbar - \hbar^{-1}) \frac{z}{1-z} \, F \otimes E
\right)  \in 
\cU_\hbar(\fg_w)^{\otimes 2} [[z]] \,.
$$
It is called the \emph{fusion operator}. It gives a certain 
canonical way to promote the R-matrix to a rational 
function on the 
maximal torus. 

We have a maximal torus 
\begin{equation}
\textup{K\"ahler torus} \times \Ct_q \subset \Aut \fgh
\label{max_tor}
\end{equation}
where 
$$
\textup{K\"ahler torus} = \Pic(X) \otimes \Ct = 
\textup{maximal torus of $\fg$} \big/ \textup{center} \,,
$$
and $\Ct_q$ acts naturally by automorphisms of 
domains of quasimaps and on $\fgh$ via the loop rotation 
automorphism \eqref{loop_rot}.  The adjoint action 
of \eqref{max_tor} on $\cU_\hbar(\fg_w)$ factors through 
a map to its maximal torus, so $\bJ$ can be promoted
to a function on \eqref{max_tor}.  Since the larger 
torus includes a shift by $q$, the equation \eqref{Jconj}
is trivially promoted to a $q$-difference equation.

\subsubsection{}

{}From the fusion operator $\bJ$, one can make 
the following dynamical operator 
\begin{equation}
\bB_w = \mathsf{m} \left( (1 \otimes \mathsf{S}_w ) 
\, \bJ_{21}^{-1} \right) \label{defB} 
\end{equation}
where 
\begin{align*}
\mathsf{m} & = \textup{multiplication} \,,\\ 
\mathsf{S}_w  & = \textup{antipode of $\cU_\hbar(\fg_w)$} \,, 
\end{align*}
and $\bJ_{21} = (12) \, \bJ$. 

Compare this with the operator 
$$
\Casimir_\alpha  
= \mathsf{m} \,  \mathbf{r}_\alpha \,, \quad 
\mathbf{r}_\alpha \in \fg_{-\alpha} \otimes \fg_\alpha \,,
$$
in which the canonical tensor $\mathbf{r}_\alpha$ is the 
$\alpha$-component of the \emph{classical} R-matrix, see 
Section 4.8 in \cite{MO}. 

This is not how the dynamical Weyl group is defined in 
\cite{EV}, but one can find a formula equivalent to \eqref{defB} for
$\fg_w=\mathfrak{sl}(2)$ there. The abstract formula 
\eqref{defB} makes sense in complete generality, for any 
$R$-matrix. 

The main result of \cite{OS} may be summarized as follows. 
Let $\nabla\subset H^2(X,\R)$ be the K\"ahler alcove which coincides
with the negative ample cone near zero. 

\begin{Theorem}\label{tgr} 
The operators $\bB_w$ define a representation of the dynamical 
groupoid of the affine K\"ahler arrangement. We have 
\begin{equation}
\bfM_\cL = \textup{const} \, \cL \,\, \dots  
\bB_{w_3} \bB_{w_2} \bB_{w_1}  \label{MBB} 
\end{equation}
where $w_1,w_2,\dots$ is the ordered set of wall crossed on the 
way from $\nabla$ to $\nabla - \cL$, and 
\begin{equation}
\bfG  = \textup{const} \, \bB_{w'_1} \bB_{w'_2} \dots \label{GBB} 
\end{equation}
where $w'_1,w'_2,\dots$ are the walls crossed on the way 
from $\nabla$ to $-\nabla$. 
\end{Theorem}

In other words, the operators 
$\bfM_\cL$ come from the lattice part of 
the dynamical groupoid, while the glue matrix $\bfG$ is the 
dynamical analog of the longest element of the \emph{finite} Weyl
group.  

The constant factor ambiguity in Theorem \ref{tgr} is of the 
same nature as the diagonal ambiguity in \eqref{polesqmult}.

\subsubsection{}

The factorizations of the operators $\bfM_\cL$ provided 
by Theorem \ref{tgr} generalize the additive 
decomposition \eqref{polesqmult} of the operators of the 
quantum differential connection.  

It is certainly natural to expect 
parallel factorizations and a 
dynamical groupoid  action for an
arbitrary symplectic resolution $X$. Note that if 
factorizations of the required kind exists, they are unique and 
determine the whole structure. 

Some speculations of 
what this might look like will be offered below. One 
should bear in mind, however, that both the quasimap 
moduli spaces \emph{and} the operators
$\bfM_\cL$, $\bfG$ depend on the choice of the 
presentation \eqref{XGIT}. This can be seen in 
example \eqref{HilbTP}. 

\subsubsection{}
The knowledge of the quantum difference equation opens the door 
 to the reconstruction of the K-theoretic DT theory of ADE fibrations and 
toric varieties in the same way as it was done in cohomology 
\cite{MOOP}. For a technical introduction to the 
subject, see e.g.\  \cite{PCMI}.

\subsubsection{}\label{s_jointly} 
As in cohomology, one can define a flat $q$-difference 
connection on the space of equivariant parameters by 
counting twisted quasimaps to $X$. For twisted quasimaps, 
the spaces $W_i$, $Q_{ij}$ from Section \ref{s_Nak} and their 
duals $\hbar^{-1} \otimes W_i^\vee$, $\hbar^{-1} \otimes Q_{i,j}^\vee$
form fixed but possibly nontrivial bundles over the domain $C$ of the 
quasimap. The $\hbar^{-1}$-weight here may also be twisted, 
that is, replaced by a line bundle on $C$. 

On general grounds, the  quantum $q$-difference connection
commutes with these shift operators. A key step in 
in the proof Theorem \ref{tgr} it to find the quantum 
Knizhnik-Zamolodchikov connection among the shift 
operators. See Section 10 of \cite{PCMI} for details. 

If we don't shift $\hbar$, that is, if we restrict ourselves to 
twisting by a maximal torus $\bA \subset \Aut(X,\omega)$,
then the resulting $q$-difference connection is a connection 
with regular singularities on the toric variety $\overline{\bA}$ 
defined by the fan of cones $\{\fC\}\subset\Lie \bA$.  Dually, the 
quantum difference connection is a connection with 
regular singularities on $\fK$, which is a toric variety defined
by the fan of ample cones. These two connections commute but
 are 
\emph{not} regular jointly --- this will be important below. 

The two connections are conjectured to be exchanged by 
symplectic duality, up-to an explicit gauge transformation.

\section{Further directions}
\label{sec:further-directions}

\subsection{Categorical dynamical groupoid}

\subsubsection{}

The operators $\bB_w$ of the dynamical groupoid 
from Theorem \ref{tgr} are normalized so that $\bB_w\big|_{z=0}=1$ 
where $0=0_X \in \fK$ is the origin in the space of K\"ahler variables
for $X$. This is natural from the curve-counting point of view. 
One can choose a more symmetric normalization, in which the 
specialization to K-theoretic 
Bezrukavnikov's groupoid from Section 
\ref{s_Bez} is the standard \eqref{limB}. In that normalization, 
one takes the factorization of $\cL$ provided by Bezrukavnikov's
theory and distributes it over the factors in \eqref{MBB}. 

\subsubsection{} 

It is natural 
to ask for a categorical lift of the dynamical groupoid. 
The dependence on the 
K\"ahler variables may be treated either formally or not
formally. In the formal treatment, we introduce an 
extra grading by the characters $H_2(X,\Z)$ of the 
K\"aher torus, otherwise nothing is changed. 

The functors $B_w$ 
from the triangle \eqref{triangDb} are 
\emph{perverse equivalences} as in the 
work of Rouquier and collaborators  \cite{CrR}, see 
\cite{ABM,Los2,Los3}. This means, in particular, 
that the category 
$$
\Db \cX_\lambda\textup{-mod} = \cC \supset \cC_1 
\supset \cC_2 \supset \dots \supset \cC_k
$$
has a finite filtration by thick subcategories such that 
the action of $B_w$ on each successive quotient 
$\cC_i/\cC_{i+1}$ is a shift of homological degree and 
equivariant grading. 

{}From \eqref{MBB} and \eqref{GBB} 
it is natural to expect the dynamical lift $B_w(z)$ to 
be the glue matrix $\bfG_{\cC_i/\cC_{i+1}}$ on each 
successive quotient, with a shift as before. This 
requires a geometric realization of the quotient 
category $\bfG_{\cC_i/\cC_{i+1}}$. In fact, if one wants to use quasimaps 
to define the glue matrix as in \eqref{defglue}, one needs
a presentation 
$$
\cC_i/\cC_{i+1} \cong \Db \Coh_{\substack{
\textup{reductive}\\ \textup{group}}} 
\begin{pmatrix}
  \textup{lci affine}\\ \textup{prequotient} 
\end{pmatrix}
\bigg/ \textup{unstable}  \,.
$$
Perhaps such presentation for each possible quotient 
may be obtained from a GIT presentation \eqref{XGIT} 
for $X$ itself. Perhaps there are ways to define the 
curve-counting operators which are more general than 
those based on quasimap spaces. In any case, it will 
require some additional data, as our running example 
\eqref{HilbTP} shows. 

\subsection{Elliptic theory}
\label{s_ell}

\subsubsection{}\label{s_ell1} 

The Yang-Baxter equation \eqref{YB} has a generalization,
in which the matrix $R(u)$ is allowed to depend on a variable 
$z$ in the maximal torus for the corresponding group. 
This is known as the \emph{dynamical} Yang-Baxter equations 
and it has very interesting solutions in elliptic functions of $u$ and
$z$. These have to do with the theory of \emph{elliptic} 
quantum groups, which was started by Felder \cite{Fel} and has 
since grown into a very rich subject. 

As anticipated in \cite{GKV} at the dawn of the subject, one gets 
elliptic quantum groups by lifting the constructions of the 
geometric representation theory to equivariant elliptic 
cohomology. 

Elliptic stable envelopes were constructed in \cite{AO}.
 Their main conceptual 
difference with both K-theoretic and categorical 
stable envelopes is the following. While stable envelopes
from Section \ref{s_st} have a \emph{piecewise constant} 
dependence on the slope $\cL\in \Pic(X) \otimes_Z \R$, 
elliptic stable envelopes depend \emph{analytically} on 
their dynamical variable 
$$
z \in \Pic(X) \otimes_\Z E \,, \quad E = \Ct/q^\Z\,.
$$
The piecewise constant dependence is recovered in the 
$q\to 0$ limit\footnote{
Using modularity of elliptic stable envelopes, one can replace 
the $q\to 0$ limit by $q\to 1$. The important part is 
that the elliptic curve $E$ degenerates to a nodal rational 
curve.} if 
$$
-\,  \Re \, \frac{\ln z}{\ln q} \to \cL  \,. 
$$
It is, of course, well-known 
that elliptic functions have a piecewise analytic asymptotics
as the elliptic curve degenerates. 

The walls across which K-theoretic stable envelopes jump 
played an absolutely key role in the preceding discussion. 
In the elliptic theory, walls become poles and all alcoves
melt into a single elliptic function. The only remaining 
discrete data is the choice of a cone $\fC \subset \Lie \bA$, and 
the choice of the ample cone in $H^2(X,\R)$ if one thinks 
about all flops of $X$ at the same time. Note that this data
is now symmetric on the two sides of the symplectic duality. 

\subsubsection{}
A categorical lift of elliptic stable envelopes along the 
lines of \eqref{stab_cat} should incorporate their 
dependence on the K\"ahler variables $z$ together with the 
elliptic modulus $q$. A few paragraphs ago we discussed 
adding such extra variables formally. A real need to 
incorporate these variables arises in certain 
categories of boundary conditions, see \cite{Kap} for a general 
introduction to the topic. 

Recall that the Nakajima varieties may be interpreted as the Higgs
branches of the moduli spaces of vacua in certain supersymmetric 
gauge theories. K-theoretic counts of 
quasimaps $C\dasharrow X$, where $C$ is a
Riemann surface, may then be 
interpreted as supersymmetric indices in 3-dimensional 
gauge theories with spacetime $C\times S^1$. 

A boundary
condition in such theory is a coupling of the bulk gauge theory
to some $(0,2)$-supersymmetric theory on $\partial C \times S^1$. 
The Hilbert space of the boundary theory is 
graded by 
$$
\textup{K\"ahler torus of $X$}  = \textup{Center of the gauge group} 
$$
together with the action of $q$ by rotations of $\partial C$. 
These Hilbert spaces  form a sheaf 
over $X$ and taking their graded index one gets an  
element of $K_\bT(X)(z)[[q]]$ which represents the $q$-expansion 
of an elliptic function.  This is how elliptic stable envelopes
should appear from a functor 
\begin{equation}
  \begin{matrix}
    \textup{boundary condition} \\
\textup{in gauge theory for $X^\bA$}
  \end{matrix}
\quad
\xrightarrow{\quad}
\quad 
 \begin{matrix}
    \textup{boundary condition} \\
\textup{in gauge theory for $X$} 
  \end{matrix}\quad ,
\end{equation}
as will be discussed in \cite{AO2}.

Note that elliptic stable envelopes depend only 
on $X$ itself and not on its presentation as a Nakajima variety. 

\subsubsection{}

While the monodromy of a differential equation gives 
a representation of the fundamental group of its regular 
locus, the monodromy of a $q$-difference equation in 
$z$ is a collection of elliptic matrices in $z$ that 
relate different fundamental solutions. In particular, 
$q$-connections with regular singularities on a toric 
variety produce elliptic connection matrices for any 
pair of fixed points. For us, those are labeled by 
a pair of birational symplectic resolutions 
$X$ for the K\"ahler
 $q$-connection or by a pair of cones in $\Lie\bA$
for the equivariant $q$-connection. 

This difference in how monodromy packages the 
information is parallel to the differences highlighted
in Section \ref{s_ell1} and, in fact, elliptic stable envelopes
give a very powerful tool for the analysis of the monodromy. 

\subsubsection{}

Recall from Section \ref{s_jointly} that the K\"ahler and 
equivariant connections are regular separately but not 
jointly. This can never happen for differential equations by 
a deep result of Deligne \cite{Del2}, but can happen very 
easily for $q$-difference equations as illustrated by 
the following system: 
$$
f(qz,a) = a f(z,a) \,, \quad f(z,qa) = z f(z,a) \,. 
$$
This proves to be a feature and not a bug, for the following 
reason. 

Let $0_X\in \fK$ be the origin in the K\"ahler moduli space for $X$
and choose a point $0_\bA \in \overline{\bA}$ on the infinity of 
the equivariant torus. Quasimap counting produces solutions 
of both K\"ahler and equivariant equations known as 
\emph{vertex functions}, or I-functions in the mirror symmetry 
vernacular. These are born as power series in $z$, the terms 
in which record the contribution of quasimaps of different 
degrees. In particular, vertex functions are 
\emph{holomorphic} in $z$ in a punctured neighborhood of $0_X$. 
However, they have \emph{poles} in $a$ accumulating to 
$0_\bA$. 

We can call these the $z$-solutions, and there are the 
dual $a$-solutions 
that are holomorphic in a punctured neighborhood of $0_\bA$. 
The connection matrix between the $z$-solutions and 
the $a$-solutions 
is an elliptic matrix in both $z$ and $a$. 
It may be called the 
\emph{pole subtraction} matrix because it may be, in principle, 
computed by canceling the unwanted poles. It may be 
compared to a monodromy matrix, or to a Stokes matrix, except 
there is no analytic continuation involved, or sectors. 
One of the main results of \cite{AO} is the identification of 
pole subtraction matrices with stable envelopes. Note 
they depend on the same discrete data, namely a choice of 
$(0_X,0_\bA) \in \fK \times \overline{A}$, and act in 
the direction 
\begin{equation}
K_\bT(X^\bA) \otimes \textup{functions}(z,a,\hbar,q)
\to K_\bT(X) \otimes \textup{same functions} \,.
\label{KXAKX}
\end{equation}
K-theories in the source and the target of \eqref{KXAKX} are 
naturally identified with initial conditions for $a$-solutions and 
$z$-solutions, respectively. 

\subsubsection{}

A combination of two steps of the form
$$
\begin{matrix}
\textup{fundamental} \\
\textup{solution at $0_\bA$} 
\end{matrix}
\quad 
\xrightarrow{\quad\textup{elliptic stable envelope}\quad} 
\quad 
\begin{matrix}
\textup{fundamental} \\
\textup{solution at $0_\bX$} 
\end{matrix}
$$
connect any two fixed points of $\fK$ or $\overline{\bA}$. 
As a result, elliptic stable envelopes: 
\begin{enumerate}
\item[---] compute the monodromy of the equivariant 
equations, 
\item[---] intertwine the monodromy of the K\"ahler 
equations for $X$ and $X^\bA$. 
\end{enumerate}
The second constraint is particularly powerful for 
Nakajima varieties, as it determines the behavior 
of the monodromy with respect to tensor product. 

\subsubsection{}

We recall that K\"ahler difference equations depend not
only on $X$ but also on a presentation of $X$ as a Nakajima 
variety. In particular, even if $X^\bA$ is a bunch of points, its 
K\"ahler difference equation may still be a
nontrivial abelian difference 
equation. As a result, even through elliptic stable envelopes
depend on $X$ only, the monodromy depends on its 
quiver presentation. 

\subsubsection{}

One of the main results of \cite{BO} is the identification of 
the categorical stable envelopes \eqref{stab_cat} with 
the parabolic induction functors for the 
quantizations under the derived equivalences 
\eqref{Dequiv}.  In particular, in the case 
$\dim X^\bA = 0$, this gives a complete control 
on the action of the triangles \eqref{triangDb} 
in K-theory and thus lets one proof the conjectures 
of Section \ref{s_Bez} in the case of isolated fixed points.

\newpage 

\noindent 
Andrei Okounkov\\
Department of Mathematics, Columbia University\\
New York, NY 10027, U.S.A.\\

\vspace{-12 pt}

\noindent 
Institute for Problems of Information Transmission\\
Bolshoy Karetny 19, Moscow 127994, Russia\\

\vspace{-12 pt}

\noindent 
Laboratory of Representation
Theory and Mathematical Physics \\
Higher School of Economics \\ 
Myasnitskaya 20, Moscow 101000, Russia 

\end{document}